\documentclass[a4paper]{article}
\usepackage[latin1]{inputenc}
\usepackage[T1]{fontenc}
\usepackage{amsmath,amssymb,stmaryrd}
\usepackage[all]{xy}
\usepackage{graphicx}

\CompileMatrices

\def\build#1_#2^#3{\mathrel{\mathop{\kern0pt#1}\limits_{#2}^{#3}}}

\def\addtrackt{<\! \mathbf{A}\!>}
\def\addtracke{<\! \mathbf{E}\!>}
\def\addtracka{<\! \mathbf{T}\!>}

\def\proof{\noindent {\it Proof.~}}

\def\fin{\hfill $\square$}

\def\part{\partial}

\def \ie {{\it i.e.~}}

\def \lra {\longrightarrow}

\def \imp {\Rightarrow}

\newcommand{\n}{\mathbb{N}}

\newcommand{\z}{\mathbb{Z}}

\numberwithin{equation}{section}

\newtheorem{theo}[subsection]{Theorem}

\newtheorem{cor}[subsection]{Corollary}

\newtheorem{lem}[subsection]{Lemma}

\newtheorem{prop}[subsection]{Proposition}

\newtheorem{defi}[subsection]{Definition}

\newtheorem{rem}[subsection]{Remark}

\newtheorem{ex}[subsection]{Example}

\def\calc{\mathcal{C}}

\def\calt{\mathcal{T}}

\def\Nil{\mathsf{Nil}}

\def\sla{\rightsquigarrow}

\def\nil{\mathsf{nil}}

\def\ab{\mathsf{ab}}

\def\gr{\mathsf{gr}}

\def\m{\medskip}

\def\trackt{\mathbf{A}}

\def\enil{\mathbf{F}_n}

\def\fnil{\mathbf{F}_m}

\def\gnil{\mathbf{F}_q}

\def\e2nil{\mathbf{F}_{2n}}

\def\group{\mathsf{Gr}}

\def\abelian{{\mathsf{Ab}}}

\def\nun{\underline{\bf n}}

\def\num{\underline{\bf m}}

\begin{document}

\title{The $\Gamma$-structure of an additive track category}
\author{G. Gaudens}
\date{August 2009}
\maketitle

\abstract{We prove that an additive track category with strong
coproducts is equivalent to the category of pseudomodels for the
algebraic theory of $\nil _2$ groups. This generalizes the classical
statement that the category of models for the algebraic theory of
abelian groups  is equivalent to the category of abelian groups. Dual statements
are also considered.}

\section{Introduction}

\bigskip  We explore in this paper the structure of \emph{additive track categories}. A track category is a $2$-category $\calc$ in which every $2$-morphism (called a track) is invertible. One can define a \emph{homotopy relation} on the morphisms of a track category: two morphisms are homotopic if there exists a $2$-morphism between them. We  obtain in this way the homotopy category $\mathrm{ho} (\calc)$ of the track category $\calc$ by identifying homotopic morphisms, and there is a canonical projection functor $\calc \lra \mathrm{ho} (\calc)$ from the underlying category of the track category to its homotopy category. Topology provides many examples of track categories: any pointed closed model category yields a track category by considering the full subcategory of fibrant and cofibrant objects, with tracks ($2$-morphisms) the homotopy classes of homotopies. Under mild conditions, the track category we obtain is part of a \emph{linear track extension}. That is, the set of self-tracks of a map depends in some functorial way on the homotopy class of the map \cite{baues1}.

In the present work, we consider \emph{additive track categories}, which are linear track extensions whose homotopy category is additive, and whose track structure is parametrized by a bilinear bifunctor on the homotopy category. Examples of additive track categories arise naturally by considering the track extension associated to a \emph{stable model category} \cite{hov}, as for example the stable homotopy category considered as the homotopy category of the Bousfield-Friedlander model category on spectra \cite{bf}. There is a natural notion of equivalence of linear track extensions, and one can  wonder if for a given additive track category, there is another one with nicer properties in its equivalence class. This question has been studied in \cite{bauesjibpira, bauespira}. In particular, it is shown in \cite{bauesjibpira} that any \emph{additive track category} has within its equivalence class an additive track category that has either strict coproducts or strict products (this terminology is explained in the appendix). We refer to this as \emph{the strictification theorem}, and any model with either strict products or strict coproducts is called a \emph{semi-strictification} in the sequel.

In this work, we explore further strictification results. Consider any object  $X$ in the homotopy category of an additive track category. This object is an Abelian group object. Can we lift this structure to the track category in some reasonable fashion? We introduce the track category $\mathbf{Pseudo}^\amalg  (\mathbf{T} , \mathbf{A})$ of pseudomodels in the track category $\trackt$ for any algebraic theory $\mathbf{T}$, which may be thought of as a model for this algebraic theory up to coherent homotopy. We derive from \cite{bauesjibpira} (more specifically \cite{p2}) that for all $n > 1$, the Abelian structure of $X$ lifts to a pseudomodel structure over the algebraic theory of $\nil _n$ groups. We also introduce the notion of a coproduct preserving pseudo natural transformation between two pseudofunctors and show that one can lift homotopy classes of maps in an essentially unique way to coproduct preserving natural transformations. There is a notion of homotopy of pseudo natural transformations of pseudofunctors and our results assemble to show that an additive track category is equivalent (as a $2$-category) to the $2$-category of pseudomodels over the category of $\nil_n$ groups. More precisely, we prove the following theorem:
\begin{theo}
\label{maintheo}
Let $\trackt$ be an additive track category with strict coproducts. For $n\geq 2$, there is a weak $\nil_n$ ringoid structure $X\lra F_X$ on $\trackt$, and the assignment:
\begin{eqnarray*}
T :&\relax \mathbf{A}\lra&\relax \mathbf{Pseudo}^\amalg  (\nil_n , \mathbf{A}) \\
&X\longmapsto &  F_X
\end{eqnarray*}
extends uniquely to a track functor $T$, which is an equivalence of track categories, and the inverse of which is the evaluation of a pseudofunctor $G :\nil _n \sla \mathbf{A}$ on the group $\mathbb{Z}$. Moreover, at the level of homotopy categories, $T$ induces the canonical equivalence
$$
\mathrm{ho} (T) : \mathrm{ho} (A) \lra \mathrm{mod}^{\amalg} (\ab, \mathrm{ho}(\mathbf{A})) .
$$
In particular, all  weak  $\nil_n$ ringoid  structures are canonically equivalent. Moreover, for $s\geq t\geq 2$, the functor $\nil_s \lra \nil_t$ induces isomorphisms
of additive track theories  $\mathbf{Pseudo}^\amalg  (\nil_t , \mathbf{A}) \lra \mathbf{Pseudo}^\amalg  (\nil_s , \mathbf{A})$.
This statement holds for $n, t=1$ under the conditions that the coefficients of $\trackt$ have no $2$-torsion. 
\end{theo}

\m The words appearing in the main theorem are explained in Section 2 and in  Appendix~A. This functor $T$ is called the \emph{$\Gamma$-structure of the
additive track category $\addtrackt$}. We do not give here any application of
the theory of $\Gamma$-structures, but rather defer it to further
papers. The $\Gamma$-structure of an additive track category is a
fundamental piece of structure and deserves therefore an independent 
treatment.

\m The paper is organized as follows. We quickly review  the notion of
algebraic theories and their models and then introduce the
$2$-categorical analogue which we call \emph{pseudomodels}
over an algebraic theory. This  allows us to state our main
results (Section \ref{additivtracketpseudo}) and prove Theorem
\ref{maintheo}, assuming certain results that are proved later. In
the next section, we show that there exists  enough pseudomodel
structures for the algebraic theories  $\nil _n$ for $n\geq 2$,
but in general not for $\nil _1=\ab$. This allows us to restate
Theorem \ref{maintheo} in terms of the notion of $\Gamma$-tracks
(Section \ref{gamma:tracks:section}), and the proof of Theorem
\ref{maintheo} is completed in Section \ref{existuniquegamma}
using this reformulation. We finish with a brief discussion of the
dual setting (Section \ref{dual}). In the appendix, we recall for the sake of self-containedness the
needed facts on track theory and pseudofunctors.

\m {\it Acknowledgements.~}This work owes a lot to M. Jibladze, to T. Pirashvili, and above all to H.-J. Baues, who shared their vast knowledge with me.
It is a pleasure to thanks them all.

\section{Additive track categories and pseudomodels}
\label{additivtracketpseudo}
\label{pseudomodsec}

\subsection{Algebraic theories}
\label{algtheo}\label{grouptheo}

The material of this section is classical, and we refer the reader to \cite{borceux} for a more complete exposition.
An algebraic theory ${\bf T}$  is a small category with a countable set of objects $\{{\bf T}_0, {\bf T}_1, \ldots ,{\bf T}_n, \ldots \} $ and specified isomorphisms from ${\bf T}_n$ to the $n$-fold categorical product of ${\bf T}_1$. In particular, ${\bf T}_0$ is a terminal object. The category with objects are \emph{product preserving} functors from ${\bf T}$ to the category of sets, and with morphisms are the natural transformations is termed the category of \emph{models for {\bf T}} and denoted by  $\mathrm{mod}^{\Pi}({\mathbf{T}})$.
We will also use the concept of a coalgebraic theory, simply defined by saying that a small category $ {\bf T}$ is a coalgebraic theory if $ {\bf T}^\mathrm{op}$ is an algebraic theory.

\m\noindent {\it Example.~}Let $\group$ be the category of groups and group homomorphisms. Let $\mathbf{F}_n$ be the free group on the set $\nun =\{1, \ldots ,n\}$. By convention, we set $\mathbf{F}_0$ to be the trivial group. We denote by $\gr$ the full subcategory of $\group$ with objects the free groups $\{\mathbf{F}_n\}_{n\in\n}$. The free group $\mathbf{F}_n$ is canonically the $n$-fold sum  of $\mathbf{F}_1 = \mathbb{Z}$ in $\group$, and the category $\gr ^{op}$ is therefore an algebraic theory, which is called the (algebraic) theory of groups.

\m\noindent {\it Example.~} Let $\abelian$ be the full subcategory of $\group$ with objects the Abelian groups. This category has products and coproducts.   Let $\mathbf{F}_n^{\ab}$ be the free Abelian group on the set $\nun =\{1, \ldots ,n\}$. In the category $\abelian$, $\mathbf{F}_n^{\ab}$ is canonically isomorphic to both the $n$-fold product of $\mathbf{F}_1^{\ab}=  \mathbb{Z}$ and the $n$-fold coproduct of $\mathbf{F}_1^{\ab}$. This means that if we let $\ab$ be the full subcategory of $\abelian $ whose objects are $\{{\bf F}_0^{\ab}, {\bf F}_1^{\ab}, \ldots , {\bf F}_n^{\ab}, \ldots \} $, then both $\ab$ and $\ab ^{\mathrm{op}}$ are algebraic theories. Reflecting the fact that finite products are isomorphic to finite coproducts of Abelian groups (within the category $\abelian$), the algebraic theories $\ab$ and $\ab ^{\mathrm{op}}$ are isomorphic.

\m\noindent {\it Example.~} A group $G$ is endowed with a natural filtration $\{\Gamma_n G\}_{n\in \n}$, \emph{the lower central series}, defined inductively by
$$
\Gamma _0 G = G, ~\Gamma _{n+1} G= [G,  \Gamma _{n} G]
$$
Here $[-,-]$ stands for commutators in $G$. A nilpotent group of class $n$ is a group $G$ such that $\Gamma _{n} G =\{0\}$. For $n\geq 1$, let  $\Nil _n$  be the full subcategory of $\group$ whose objects are  groups of nilpotency $n$. The case of $n=1$ corresponds to the category  of Abelian groups, and the category $\Nil _1$ is simply denoted by $\abelian$, while the case $n=2$ corresponds to the case of so called $nil$-groups which we denote simply by $\Nil$. The categories $\Nil _n$ have free objects on arbitrary sets of generators. Indeed, for any set $S$, denote the free group on $S$ by $<\!S\!>$. The group $<\!S\!>_n=<\!S\!>/\Gamma _{n} <\!S\!>$ is a nilpotent group of class $n$ such that the adjunction formula
$$
\mathit{Sets} (S, G)\cong \Nil  _n (<\!S\!>_n, G)
$$
holds for all sets $S$ and nilpotent group $G$ of class $n$. In other words, $<\!S\!>_n$ is the free object generated by $S$.

Let $\nil  _n$ be the full subcategory of $ \Nil  _n$ with objects the  countable set $\{\mathbf{F}^{\nil  _n}_{p}\}_{n \in \n}$, where $ \mathbf{F}^{\nil  _n}_{p}$ is the free nilpotent group of class $n$ on the finite set $\{1, \ldots, p\}$ and $ \mathbf{F}^{\nil  _n}_{0}$ is the trivial group by convention. For $n=1$, we recover $\nil _1 = \ab$. The group $ \mathbf{F}^{\nil  _n}_{p}$ is canonically isomorphic to the $n$-fold coproduct of $ \mathbf{F}^{\nil  _n}_{0}$ in $\nil  _n$, so that the category $\nil  _n^{op}$ is an algebraic theory.

\m These examples relate in the following way. The lower central series induces an augmented tower of coproduct preserving functors
$$
\group \lra  \{\ldots\lra  \Nil (n+1) \lra  \Nil  _n \lra \ldots \lra  \abelian\}
$$
where the functor $\group \lra \Nil  _n$ maps a group $G$ to $G/\Gamma _{n+1} G$. That is to say, we get an augmented tower of algebraic theories
\begin{eqnarray}
\label{tower}
\gr^{op} \lra  \{ \ldots \lra \nil _{n+1}^{op} \lra  \nil  _n ^{op}\lra \ldots \lra \ab ^{op} \} \quad .
\end{eqnarray}
In the sequel, we consider $\gr$ as $\nil_n$ with $n=+\infty$.

Let $\calc$ be a category with finite products and ${\bf T}$ an algebraic theory. We let $\mathrm{mod}^{\Pi} ({\bf T}, \calc)$ be the category of product preserving functors ${\bf T}\lra \calc$ and their natural transformations, termed the category of models in $\calc$ for the algebraic theory ${\bf T}$.
We will also consider the category $\mathrm{mod}^{\amalg} ({\bf T}^{\mathrm{op}}, \calc)$ of coproduct preserving functors ${\bf T}^{\mathrm{op}}\lra \calc$ and their natural transformations, termed the category of coproduct models (or $\amalg$-models) in $\calc$ for the coalgebraic theory ${\bf T}^{\mathrm{op}}$. If $\calc$ is additive, there are equivalences of categories
$$
\mathrm{mod}^{\amalg} (\ab^{\mathrm{op}}, \calc)\cong \mathrm{mod}^{\Pi} (\ab^{\mathrm{op}}, \calc)\cong \mathrm{mod}^{\amalg} (\ab, \calc)\cong \mathrm{mod}^{\Pi} (\ab, \calc)
$$
and we recall that:
\begin{prop}
\label{carac:cat:add}
A category $\calc$ is additive if and only if the evaluation functor
$$
\mathrm{mod}^{\Pi} (\ab^{\mathrm{op}}, \calc) \lra \calc
$$
is an equivalence of categories.
\end{prop}

\subsection{Pseudomodels for additive track categories}
\label{pseudomod}

We now consider the case of an additive track category $\trackt$ (see the appendix for the notations):
$$
\xymatrix{
D+~\ar@{>->}[r]&\relax \mathbf{A}_1\ar@<0.3ex>[r]\ar@<-0.3ex>[r]&\relax \mathbf{A}_0 \ar@{->>}[r]^{p}&\relax \mathrm{ho}(\mathbf{A})= \calc .\\
}
$$
The homotopy category
of $\trackt$ is additive, therefore every object in $\trackt$ has the structure of a group and a cogroup \emph{up to homotopy}.
It is a result of \cite{bauesjibpira} (see the appendix) that any additive track category has in its equivalence class an additive track category having either strict products or strict coproducts. There is nevertheless an obstruction to obtain both properties simultaneously, as we shall see. The main point here is that from a theoretical viewpoint, there is no loss of generality in considering track categories having either strict products or strict coproducts.

\m\noindent{\bf Convention.~~} {\it In the sequel, we will consider the case of \emph{additive track categories with strict coproducts} in order to fix the ideas, except otherwise stated (for instance in Section \ref{dual} where we state the dual  results).}

\subsubsection{Coproduct preserving pseudofunctors and natural transformations}

The necessary material on pseudofunctors is recalled in the appendix (see \ref{twocatpseudo}).
Let $\calc$ and $\calt$ be  additive track categories with strict coproducts.
A pseudofunctor $\varphi :\calc \sla \calt$ is coproduct preserving if the following conditions are satisfied:

\begin{itemize}
\item $\varphi$ is both reduced at tracks and objects (\ie \emph{completely reduced}, see \ref{pseudopseudo}),
\item the natural map $\iota_{X, Y}= \varphi (X\lra X\vee Y)\vee \varphi (Y\lra X\vee Y) $
\begin{eqnarray}
\iota_{X, Y}: \varphi(X)\vee \varphi (Y) \lra \varphi (X \vee Y)
\end{eqnarray}
is an isomorphism for all $X$ and $Y$ in $\calc$,
\item for all objects $X, Y, Z, T$ and maps $h: X\lra Z, k: Y\lra Z, g:Z\lra T$ in $\calc$, the equation
\begin{eqnarray}
\varphi_{g,(h\vee k)} \iota _{X, Y} = \varphi_{g,h}\vee \varphi_{g,k}
\end{eqnarray}
is satisfied,
\item for all objects $X$ in $\calc$, the equation
\begin{eqnarray}
\varphi_{X\vee X}= (\varphi_X , \varphi_X) \quad .
\end{eqnarray}
\end{itemize}

Let $\varphi, \psi: \calc \sla \calt $ be two coproduct preserving pseudofunctors. We say that a pseudo natural transformation $\alpha : \varphi \lra \psi$  is coproduct preserving if the following conditions are satisfied:
\begin{itemize}
\item for all $X, Y$, the following diagram is commutative
\begin{eqnarray}
\xymatrix{
\relax \varphi (X)\vee \varphi (Y)\ar[r]^{\iota_{X, Y}}\ar[d]_{\alpha _X \vee \alpha_Y}&\varphi (X\vee Y)\ar[d]^{\alpha _{X\vee Y}}\\
\relax \psi (X)\vee \psi (Y)\ar[r]^{\iota_{X, Y}}&\psi (X\vee Y) &,\\
}
\end{eqnarray}
where $\iota_{X, Y}$ is as above and moreover,

\item for all $X, Y, Z$ and all maps $f:X\lra Z, G: X\lra Z$, the pasting of tracks in the diagram (\ref{diagdiag}) yields $\alpha _{f}\vee \alpha_{g}$.
\end{itemize}

\begin{eqnarray}
\label{diagdiag}
\xymatrix{
\relax \varphi (X)\vee \varphi (Y)\ar[r]^{\iota_{X, Y}}\ar[d]_{(\alpha _X ,\alpha_Y)}&\ar@{}[rd]|{\Downarrow \alpha _{f\vee g}}\varphi (X\vee Y)\ar[d]_{\alpha _{X\vee Y}}\ar[r]^{\varphi (f \vee g)} & \varphi (Z) \ar[d]^{\alpha _Z} \\
\relax \psi (X)\vee \psi (Y)\ar[r]_{\iota _{X, Y}}&\psi (X\vee Y)\ar[r]_{\varphi (f \vee g)} & \psi (Z)\\
}
\end{eqnarray}
Now given two coproduct preserving maps $\alpha, \beta : \varphi \lra \psi$ of coproduct preserving pseudofunctors, a homotopy is called coproduct preserving if the pasting in the following diagram yields the identity track:
\begin{eqnarray}
\xymatrix{
\ar@{}[rd]|{\relax \build{}_{H}^{\Rightarrow}}&\relax \varphi (X)\vee \varphi (Y)\ar[r]^{\iota_{X, Y}}\ar[d]^{(\alpha _X , \alpha_Y)}\ar@/_3pc/[d]_{(\beta _X , \beta_Y)}
&\varphi (X\vee Y)\ar[d]_{\alpha _{X\vee Y}}
\ar@/^3pc/[d]^{\beta _{X \vee Y}}_{\relax \build{\Rightarrow}_{~~~ H_{X\vee Y}}^{}}
\\
&\relax \psi (X)\vee \psi (Y)\ar[r]_{\iota_{X, Y}}&\psi (X\vee Y)&,\\
}
\end{eqnarray}
where $H$ is the track $(H_X^\boxminus , H_Y^\boxminus )$.
The following proposition is a direct consequence of the definitions.
\begin{prop}
Let $\calc$ and $\calt$ be two track categories with strict coproducts. Assume that $\calc$ is a small track category. Then coproduct preserving pseudofunctors, coproduct preserving pseudo natural transformations and coproduct preserving homotopies build a track subcategory $\mathbf{Pseudo} ^\amalg(\calc, \calt)$  of   $\mathbf{Pseudo} (\calc, \calt)$.
\end{prop}
Recall that the superscript $\boxminus$ indicates the inverse of a track.

\subsubsection{Coproduct pseudomodels}
\label{statements}

Let us assume that $\trackt$ is  a track category with strict coproducts. We consider a coalgebraic theory $\mathbf{T}$, considered as a discrete track category (that is, with only identity tracks).
\begin{defi}
The category of $\amalg$-pseudomodels in $\trackt$ for the algebraic theory $\mathbf{T}$ is the category $\mathbf{Pseudo} ^\amalg(\mathbf{T}, \trackt)$.
\end{defi}
We first notice that in this setting the notion of  coproduct preserving homotopy of coproduct preserving pseudofunctors is particularly simple:
\begin{prop}
\label{so:simple}
Let $\calt$ be a track category with strict coproducts and $\mathbf{T}$ be an algebraic theory. Given two coproduct preserving pseudofunctors
$F, G: \mathbf{T} \sla \calc$.
The set of coproduct preserving homotopies $F\imp G$ is the set of tracks $F(\mathbf{T}_1)\imp G(\mathbf{T}_1)$, obtained {\it via} the evaluation functor.
\end{prop}

The proof is staightforward and will therefore be omitted.
In the following, we are interested in the case $\mathbf{T} = \nil _n$, for which we introduce further terminology.
\begin{defi}
If $\trackt$ is a track category with strict coproducts, then for $n\geq 1$,
\begin{itemize}
\item A \emph{weak $\nil_n$ cogroup structure} on the object $X$ of $\trackt$ is a  $\amalg$-pseudomodel $\varphi$ in $\trackt$ for the algebraic theory $\nil _n$, such that $\varphi (\mathbb{Z})= X$. A \emph{weak $\nil _n$ cogroup} in $\trackt$ is a couple $(X, F_X)$ where $X$ lies in $\calc$ and $F_X$ is a weak  $\nil _n$ cogroup structure on $X$. A map of weak  $\nil _n$ cogroups $(f, \Gamma ^f):X \lra Y$ is a map $f:X\lra Y$ in $\trackt$ together with a coproduct preserving pseudo natural transformation $\Gamma ^f : F_X \lra F_Y$ such that the evaluation
$$
\Gamma ^f (\mathbb{Z}) : X= F_X (\mathbb{Z})\lra F_Y (\mathbb{Z}) =Y
$$
is equal to $f$.
\item A \emph{semi $\nil_n$ ringoid structure} on $\trackt$ consists of a chosen weak  $\nil _n$ cogroup structure for every object $X$ in $\trackt$.
\end{itemize}
\end{defi}
We first notice:
\begin{prop}
Let $\trackt$ be a track category with strict coproducts whose homotopy category $\mathrm{ho}(\trackt)$ is additive, let $(X, F_X)$ be a weak  $\nil _n$ cogroup in $\trackt$. Then $X$ projects to the Abelian cogroup structure in the homotopy category of $\trackt$.
\end{prop}
That is, let $F_X : \nil_n \lra X$ be a weak  $\nil _n$ cogroup structure on $X$. This induces a weak  $\nil _n$ cogroup structure on $X$ viewed as an object of $\mathrm{ho}(\trackt)$. This means in particular that  $X$ is a cogroup object in two different ways, one coming from the additive structure of $\mathrm{ho}(\trackt)$, and
one coming from $F_X$. Moreover, $F_X$ is a cogroup structure \emph{in the additive category $\trackt$}. A well known trick shows that in such a situation, both structure have to coincide. That is to say, $pF_X$ factors through $\ab$ (where $p: \trackt \rightarrow \mathrm{ho}(\trackt)$ is the projection functor) and this factorization is precisely the Abelian structure $\varphi _X$ of $X$ in $\mathrm{ho}(\trackt)$. In other words, there is a commutative diagram
$$
\xymatrix{
\relax \nil_n \ar[r]^{\relax F_X}\ar[d]&\relax\trackt \ar[d]^{p}\\
\relax \ab \ar[r]_{\varphi_X}&\relax \mathrm{ho}(\trackt)&.
}
$$
Assume now that $\trackt$ is the underlying track category of a semi-strictified linear track extension $\addtrackt$. A weak $\nil_n$ cogroup structure on $X$ descends in particular to
a cogroup structure on $X$ in the homotopy category. But we assume that $\addtrackt$ has an additive homotopy category, hence there is a factorization of
$F_X $ through $\ab$, and according to proposition \ref{carac:cat:add}, the Abelian cogroup structure on $X$ is unique up to isomorphism, thus \emph{any weak $\nil _n$ cogroup structure on $X$ lifts the natural structure in the homotopy category}.
It is therefore  natural to ask whether an arbitrary object in an additive track category with strict coproducts,
which is automatically a cogroup and a group in the homotopy category,  is a cogroup in $\trackt$ or not. The answer is:
\begin{theo}
\label{exist}
Any object $X$ in an additive track category $\addtrackt$ with strict coproducts admits the structure of a weak $\nil_n$ cogroup $F_X$, for $2\leq n\leq +\infty$. This statement holds for $n=1$ under the extra assumption that the linear system has no $2$-torsion.
\end{theo}
We notice that any track category with Abelian $2$-automorphism groups (morphisms form Abelian groupoids) is canonically and in an essentially unique way part of a linear track extension, hence Theorem  \ref{exist}
applies as well to such categories provided that the homotopy category is additive (\cite{bauesjib}, see also Section \ref{strictproducts} of the appendix), and the associated natural system is a bilinear bifunctor. Hence for every object $X$ there is a weak $\nil_n$ cogroup structure $F_X$ on $X$, for $n\geq 2$,
$$
\xymatrix{
\relax \nil_n \ar@{.>}[r]^{F_X}\ar[d]&\relax\trackt \ar[d]\\
\relax \ab \ar[r]_{\varphi_X}&\relax \mathrm{ho}(\trackt)
}
$$
where $\varphi_X$ is the Abelian structure coming from the additive structure of $\mathrm{ho}(\trackt)$.
We notice that we have passed from $\ab= \nil_1$ to $\nil_n$ with $n\geq 2$ to get such a lifting (this is explained in Section \ref{existenceofpseudo}).
\begin{cor}
Any additive track category with strict coproducts can be given the structure of a weak $\nil_n$ ringoid.
\end{cor}
We provide a proof of this fact in  Section  \ref{existenceofpseudo}. We then prove (Section \ref{existuniquegamma}):
\begin{theo}
\label{uniqnat}
Let $X, Y$ be two objects in a  weak $\nil_n$ ringoid $\addtrackt$, having the weak $\nil_n$ cogroup structures $F_X$ and $F_Y$. Then any map $f:X\lra Y$ extends uniquely to a map
$$
\Gamma ^f: F_X \lra F_ Y
$$
in the category $\mathbf{Pseudo}^{\amalg} (\gr , \mathbf{A})$.
\end{theo}
These results assemble together to produce a proof of Theorem \ref{maintheo}

\m {\it Proof of Theorem \ref{maintheo}.~} We will prove  Theorem \ref{maintheo} assuming Theorem \ref{exist} and \ref{uniqnat}. We first notice that $G$ is clearly a track functor. Theorem \ref{exist} and \ref{uniqnat} say that  $T$ extends uniquely to a track functor
$$T : \mathbf{A}\lra \mathbf{Pseudo}^\amalg (\nil _n , \mathbf{A}) ~~.$$
Any track $s: f \imp g$ in $A$ can be extended to a coproduct preserving homotopy $\Gamma ^s: \Gamma^f \imp \Gamma^g$, and this in a unique way (by proposition \ref{so:simple}). We therefore get a track functor
$$
T : \mathbf{A}\lra \mathbf{Pseudo}^\amalg (\nil _n , \mathbf{A})\quad ,
$$
such that the composition $GT: \mathbf{A} \lra \mathbf{A}$ is the identity functor.
On the other hand, given an object $F$ in $\mathbf{Pseudo}^\amalg (\nil _n , \mathbf{A}) $, the object $a=G(F)$ in $\mathbf{A}$ comes with a weak $\nil_n$ cogroup structure defined by $F$, which might or might not coincide with the chosen one $F_a$ (that of Theorem \ref{existgamma}). Nevertheless, by Theorem \ref{existgamma}, the identity $a\lra a$ extends uniquely to a pseudo natural transformation $F \sla F_a$, and this transformation is an isomorphism. This finishes the proof of the theorem.

\fin

\section{Existence of pseudomodels}
\label{existenceofpseudo}

The aim of this section is to construct a weak $\nil_n$ cogroup structure for every object in an additive track category with strict coproducts, that is to prove Theorem \ref{exist}. In other words, we show that any additive track category with strict coproducts is a  weak $\nil_n$ ringoid. The proof is  based on cohomological arguments. We mention that there is an alternative  direct constructive proof of  Theorem
\ref{exist}, which is in the spirit of the proofs in Section \ref{existuniquegamma} (see \cite{groupcase}). This approach has only been amenable to the author for $n=+\infty$, and will therefore
be omitted here.  Let us recall some facts on the \emph{cohomology of categories} (\cite{bauesd, baueswir}, see also the appendix).

\subsection{Cocycles for linear track extensions}
Recall that all the necessary notations and definitions  are given in the appendix.

Let $\addtracke$ be a linear track extension of the small category $\calc$ by the natural system $D$ (see Section  \ref{lineartrack}), and let $\mathbf{E}$ be the underlying track category, $\mathbf{E}_0$ and $\mathbf{E}_1$ as in \ref{trackcatpar}
$$
\xymatrix{
D+~\ar@{>->}[r]& \relax \mathbf{E}_1\ar@<0.3ex>[r]\ar@<-0.3ex>[r]&\relax \mathbf{E}_0 \ar@{->>}[r]^{p}&\relax \mathrm{ho}(\mathbf{E})= \calc\\
}
$$
We choose  functions:
$$
t: \mathrm{Mor}(\calc) \lra \mathrm{Mor}(\mathbf{E}_0) ,~~~ H: N_2(\calc) \lra \bigcup_{f, g \in Mor (\mathbf{E}_0)}  \llbracket f, g \rrbracket
$$
such that
\begin{itemize}
\item $t$ sends identities to identities,
\item $pt (f) =f$ for any morphism $f$ of $\calc$, and
\item $H(f, g)\in D(tf\circ tg, t(f\circ g))$.
\end{itemize}

We define a function $c_{T} (t, H): N_3 (\calc) \lra \bigcup_{(f, g, h)\in  N_3(\calc) } D(fgh)$
by
$$
c_{T} (t, H)(f, g, h)= \sigma ^{-1}_{t(fgh)}(\Delta _T (f, g, h))
$$
where $\sigma$ is the structure isomorphism of the linear track extension, see (\ref{lineartrackpar}). We also define:
$$
\Delta _T (f, g, h)= -H(f, gh)-(tf)_* H(g, h)+ (th)^* H (f, g)+ H(fg, h).
$$
\begin{lem}[B-D lemma A.1]
Let $c_{T}$ be defined as above. We have:
\begin{enumerate}
\label{lemme:bauesD}
\item $c_T (t, H)$ is a cocycle in $C^3 (\calc, D)$,
\item if $c$ is a $2$-cochain in $C^2 (\calc, D)$, then
$$
\delta c + c_T (t, H) = c_T (t, H-c),
$$
where $(H-c) (f, g)= H(f, g)- \sigma _{t(fg)} c(f, g)$,
\item The class of  $c_T (t, H)$ in $H^3 (\calc, D)$ does not depend on $t$ and $H$,
\item The class of  $c_T (t, H)$ in $H^3 (\calc, D)$ depends only on the component of $T$ in $\mathbf{Track}(\calc , D)$.
\end{enumerate}
\end{lem}

Assume $\calc$ is a small category, $D$ a natural system on $\calc$, and consider
a linear track extension $\addtracka$ of $\calc$ by $D$ such that the corresponding class in $H^3 (\calc, D)$
is trivial. By Lemma  \ref{lemme:bauesD}, as $t$ and $H$ vary, $c_T (t, H)$ describes exactly all
the possible cocycles representing the class of $\addtracka$ in $H^3 (\calc, D)$.  Hence we have, for some $t$ and $H$:
$$
c_{\mathbf{T}} (t, H)(f, g, h)= \sigma ^{-1}_{t(fgh)}(\Delta _{\mathbf{T}} (f, g, h))
$$
and
$$
0= \Delta _{\mathbf{T}} (f, g, h)= -H(f, gh)-(tf)_* H(g, h)+ (th)^* H (f, g)+ H(fg, h).
$$
This equality means that $(t,H): \calc\lra \mathcal{D}$ is a \emph{pseudofunctor}. That is:
\begin{cor}
\label{existsect}
Let $\addtracke$ be a linear track extension
$$
\xymatrix{
D+~\ar@{>->}[r]&\relax \mathbf{E}_1\ar@<0.3ex>[r]\ar@<-0.3ex>[r]&\relax \mathbf{E}_0 \ar@{->>}[r]^{p}&\relax \mathrm{ho}(\mathbf{E})= \calc.\\
}
$$
Then the associated class in $H^3 (\calc, D)$ is trivial if and only if there is a pseudofunctor $(t,H):\calc\lra \mathbf{E}$ such that $pt = id$.
Moreover, if $\addtracke$ is semi-strictified, then  $t$ is a \emph{coproduct preserving} pseudofunctor. Such a pseudofunctor is called a \emph{pseudosection} in
the sequel.
\end{cor}

\subsection{Construction of pseudomodels}

Let $\addtracke$ be a semi-strictified track extension of the small category $\calc$ by the natural system $D$ (see Section  \ref{lineartrack}), and let $\mathbf{E}$ be the underlying track category, $\mathbf{E}_0$ and $\mathbf{E}_1$ as in Section \ref{trackcatpar}
$$
\xymatrix{
D+~\ar@{>->}[r]&\relax \mathbf{E}_1\ar@<0.3ex>[r]\ar@<-0.3ex>[r]&\relax \mathbf{E}_0 \ar@{->>}[r]^{p}&\relax\mathrm{ho}(\mathbf{E})= \calc .\\
}
$$
Let $c$ be an object in $\calc$ and $\varphi_c$ be the Abelian cogroup structure on $c$, and $\nil _n$ be any of the
algebraic theories described in Section \ref{algtheo}, and $\lambda _n$ be the abelianization functor. We have a commutative diagram
\begin{eqnarray}
\label{wearethechamp}
\xymatrix{
\relax \lambda_n ^*\varphi_c^* X \ar[r] \ar[d] &\varphi_c^* X\ar[r]\ar[d]&\relax \mathbf{E}\ar[d]\\
\ar@{.>}@/_7pt/[u]_{(i)}\ar@{.>}[ru]_{(ii)}\relax \nil _n \ar[r]_{\lambda _n} &\relax\ab\ar[r]_{\varphi_c} &\relax\mathrm{ho} (\mathbf{E})= \calc &.
}
\end{eqnarray}
Because this diagram is a pullback diagram of categories, the existence of a pseudosection $(i)$ or the existence of the lifting 
$(ii)$ are equivalent, and moreover,  $(i)$ is coproduct preserving if and only if $(ii)$ is.
Hence, according to corollary \ref{existsect}, such a pseudosection exists if and only if the  characteristic class $\addtracke \in H^3 (\calc, D)$
is mapped to zero under the composition of natural maps
$$
H^3 (X, D) \longrightarrow H^3(\ab, \varphi_c ^* D) \longrightarrow H^3(\nil _n, p^*\varphi_c ^* D)~~.
$$
But according to \cite[Theorem A.1]{p2}  (special case $L=1$ of this theorem actually):
\begin{theo}
\label{vanish}
$H^3 (\nil _2 , D)$ is trivial for any biadditive bifunctorial natural system. If $D$ has no $2$-torsion, then $H^3 (\nil _1 , D)= H^3 (\ab, D) $ is already trivial.
\end{theo}
Hence Theorem \ref{exist} is proved, as a direct consequence of Theorem \ref{vanish}.

\section{\relax $\Gamma$-tracks}
\label{gamma:tracks:section}

In this section $\addtrackt$ is an additive track category with strict coproducts, together with a weak $\nil_s$ ringoid structure for some $s\geq 1$ (see Section \ref{existenceofpseudo} for the existence of such a structure).
After having introduced some terminology, we state in this section three theorems (\ref{existgamma}, \ref{gammarespmap} and \ref{gammaresptrack}) that altogether
complete the proof of Theorem \ref{maintheo}. Theorem \ref{gammarespmap} and \ref{gammaresptrack} are consequences of  Theorem \ref{existgamma},   the proof of which is postponed to Section \ref{existuniquegamma}.

\subsection{Some notations}

We single out certain maps in $\nil _n$  which play an important role in the following. Let $\nun\subset \num$ be an injection (an inclusion by abuse of notation). The functorially associated  map of free $\nil_s$ groups $\mathbf{F}_n^{\nil _s}\lra\mathbf{F}_m^{\nil _s}$ is called an \emph{inclusion}. An inclusion $\nun\subset \num$ defines also a map in the reverse direction $\mathbf{F}_m^{\nil _s}\lra  \mathbf{F}_n^{\nil _s}$ by taking a generator $f\in \num$ to itself if $f\in \nun$ and to $1$ otherwise. Such maps are called \emph{projections}.
Let $\alpha _{\! n} : \z \lra \enil$ be the map that sends the generator of $\z$ to the product of the generators of $\enil^{\nil _s}$ in the increasing order. For all $e\in \nun$, let $r_e:\enil^{\nil _s} \lra \z$ be the retraction on the $e^{th}$ summand.
Dually, let $\beta _{\! n} :\enil \lra  \z$ be the map that sends all generators of $\enil^{\nil _s}$ to $1\in \z$.  For all $e\in \nun$, let $i_e:\z\lra \enil^{\nil _s}$ be the inclusion of the $e^{th}$ summand in $\enil^{\nil _s}$.
a weak $\nil_n$ cogroups $(X, F)$, then the object $F(\mathbf{F}_n)$ is isomorphic to $\vee _n X$. This isomorphism is made implicit in the following as it plays no significant role. For $\alpha :\mathbf{F}_n^{\nil _s} \lra  \mathbf{F}_m^{\nil _s}$, the associated map $F(\alpha)$ is simply denoted by $\alpha$.

In the following, given objects $X, Y$ and a map $f: X\lra Y$, the notation $(f)_n: \vee_{i=1}^{n} X  \lra  \vee_{i=1}^{n} X$ means $\vee_{i=1}^{n} f$.

\subsection{ \relax $\Gamma$-structures}

Let $\addtrackt $ be an additive track category with strict coproducts.
We introduce the convenient concept of $\Gamma$-tracks, which is the local counterpart of a coproduct preserving natural transformation. We construct such $\Gamma$-tracks on $\trackt$ in a natural fashion. This leads (finally) to the proof of Theorem \ref{uniqnat}.
We first need to introduce interchange tracks.

\begin{defi}
Let $\addtrackt$ be an additive track category with strict coproducts, and $f:X\lra Y$ be a map between two weak $\nil_s$ cogroups. For $\alpha: \nun \lra \num$ in ${\nil _s}$ and $f:X\lra Y$ an interchange track is a track $\alpha(\vee _{\! n} f) \imp (\vee _{\! m} f)\alpha $ as in diagram (\ref{diag25}). An interchange structure for $f:X\lra Y$ is a correspondence
$$\alpha\longmapsto  \Gamma _{\alpha}^{f}$$
such that $\Gamma _{\alpha}^{f}$ is the trivial track as soon as $\alpha$ is either an inclusion or a projection,
\begin{eqnarray}\label{diag25}
\xymatrix{
\relax \vee _{\! n} X \ar[r]^{\alpha}\ar[d]|{\relax (f)_n} \ar@{}[rd]|{\Downarrow \Gamma _{\alpha}^{f}}&\vee _{\! m} X \ar[d] |{\relax (f)_m} \\
\relax \vee _{\! n} Y \ar[r]_{\alpha} &\vee _{\! m} Y&.
}
\end{eqnarray}
\end{defi}
To make sense of this definition, one notices that this diagram is commutative in the homotopy category, and therefore there is at least one such a track.
Suppose we are given some interchange structure $\Gamma^f$ for a map $f$ between weak $\nil_n$ cogroups.
We can introduce a new operation $\boxbox$ on $\trackt$ by pasting interchange tracks along the weak $\nil_n$ cogroup structures.
Let $\alpha:\enil^{\nil _s} \lra \fnil^{\nil _s} , \beta: \fnil^{\nil _s} \lra \gnil^{\nil _s}$ be group homomorphisms and $f:X\lra Y$ be a map in $\trackt$. We define $\Gamma ^f _{\beta} \boxbox\Gamma ^f _{\alpha}  $ to be the pasting of tracks in the following diagram:
\begin{eqnarray}
\label{diag26}
\xymatrix{\relax &&\\
\relax \vee _{\relax \!   n} X\ar@/^2pc/[rr]^{\relax\beta\alpha} \ar[r]^{\alpha}\ar[d]|{(f)_n }
\ar@{}[rru]|{\Downarrow \varphi_{\beta, \alpha}^{\boxminus}}
\ar@{}[dr]|{\relax \Downarrow \Gamma ^f _{\alpha}}
& \vee_{\relax \!   m} X  \ar[r]^{\beta} \ar[d]|{(f)_m }
\ar@{}[dr]|{\relax \Downarrow \Gamma ^{f}_{\beta}}
&\vee_{\relax \!   p} X  \ar[d]|{(f)_p}\\
\relax \vee_{\relax \!   n} Y  \ar@{}[rrd]|{\Downarrow \varphi_{\beta, \alpha}}\ar@/_2pc/[rr]_{\beta\alpha}\ar[r]^{\alpha}& \vee_{\relax \!   m} Y  \ar[r]^{\beta} &\vee_{\relax \!   p} Y   &.\\
\relax &&
}
\end{eqnarray}
Recall  that  $s^{\boxminus}$ denotes the inverse of the track $s$. 
The following proposition is an elementary consequence of the pseudofunctor property that defines a weak $\nil_s$ cogroup.
\begin{prop}
Let $\addtrackt$ be track category and let $f$ be a map between two weak $\nil_n$ cogroups in $\addtrackt$.  For any associated interchange structure, the operation $\boxbox$ on interchange tracks for $f$ is associative.
\end{prop}
\begin{defi}
Let $\addtrackt$ be an additive track category with strict coproducts.
We say that an interchange structure (associated to some map $f:X\lra Y$ of weak $\nil_n$ cogroups in $\trackt$) satisfies \emph{property $(\Gamma)$} if for all maps in $\gr$ $\alpha:\enil^{\nil _s} \lra \fnil^{\nil _s}$ and $\beta: \fnil^{\nil _s} \lra \gnil^{\nil _s}$ maps in $\nil _s$ in the additive track category $\addtrackt$ with strict coproducts, the pasting in the diagram
\begin{eqnarray}
\label{diag100}
\xymatrix{\relax &&\\
\relax \vee _{\relax \!   n} X\ar@/^2pc/[rr]^{\relax\beta\alpha} \ar[r]^{\alpha}\ar[d]|{(f)_n }
\ar@{}[rru]|{\Downarrow \varphi_{\beta, \alpha}^{\boxminus}}
\ar@{}[dr]|{\relax \Downarrow \Gamma ^f _{\alpha}}
& \vee_{\relax \!   m} X  \ar[r]^{\beta} \ar[d]|{(f)_m }
\ar@{}[dr]|{\relax \Downarrow \Gamma ^{f}_{\beta}}
&\vee_{\relax \!   p} X  \ar[d]|{(f)_p}\\
\relax \vee_{\relax \!   n} Y  \ar@{}[rrd]|{\Downarrow \varphi_{\beta, \alpha}}\ar@/_2pc/[rr]_{\relax\beta\alpha} \ar[r]^{\alpha}& \vee_{\relax \!   m} Y  \ar[r]^{\beta} &\vee_{\relax \!   p} Y   \\
\relax &&\\
}
\end{eqnarray}
yields $\Gamma_{\beta\alpha}^f$. That is, if
$$
\Gamma^f_{\beta} \boxbox \Gamma^f_{\alpha}= \Gamma^f_{\beta \alpha}
$$
for all $\alpha , \beta$. An interchange structure satisfying property $(\Gamma)$ is termed a $\Gamma$-structure associated to $f$ and the interchange tracks are termed \emph{$\Gamma$-tracks} associated to $f$.
\end{defi}
A direct consequence of the definitions is
\begin{prop}
\label{propcompsum}
Let $\addtrackt$ be an additive track category with strict coproducts and let $f$ be a map between weak $\nil_n$ cogroups. The following formula holds  for any interchange structure $\alpha\longmapsto  \Gamma _{\alpha}^{f}$ which is a $\Gamma$-structure:
\begin{eqnarray}
\Gamma ^f_{\vee _i \alpha _i}= \vee _i \Gamma ^f_{\alpha _i} \quad .
\end{eqnarray}
\end{prop}
The following remark shows the interest of the notion of a $\Gamma$-structure associated with a map of weak $\nil_s$ cogroups:
\begin{rem}
\label{remremrem}
Let $f$ be $X\lra Y$ be a map between the underlying spaces of the weak $\nil_s$ cogroups $(X, F_X)$ and $(Y, F_Y)$. Then a $\Gamma$-structure  $\Gamma ^f$ associated with $f$ is nothing but a coproduct preserving pseudo natural transformation $\Gamma ^f :F_X \lra F_Y$ between the coproduct preserving pseudofunctors $F_X, F_Y : \gr \lra \trackt$ such that $\Gamma ^f(\mathbb{Z}) = f : X= F_X (\mathbb{Z}) \lra F_Y (\mathbb{Z})= Y$.
\end{rem}

\subsection{Existence of $\Gamma$-structures}

We have the following local existence theorems.
\begin{theo}
\label{existgamma}
Let $\addtrackt$ be an additive track category with strict coproducts. Any map $f$ between two weak $\nil_s$ cogroups in $\trackt$ has a unique associated  $\Gamma$-structure.
\end{theo}
The uniqueness statement in Theorem \ref{existgamma} yields also the following two theorems \ref{gammarespmap}  and \ref{gammaresptrack}. Indeed, one simply checks that pasting yields exactly an interchange structure satisfying property $(\Gamma)$, and the equality follows from the uniqueness.
\begin{theo}[\emph{Naturality with respect to maps in $\trackt$}]
\label{gammarespmap}
Let  $\alpha:\enil \lra \fnil$ be a map in $\gr$, and  $f:X\lra Y$, $g:Y\lra Z$ be maps of weak $\nil_s$ cogroups in $\trackt$. The unique $\Gamma$-structures associated to $f$ and $g$ satisfy naturality with respect to maps in $\trackt$: the pasting in the diagram
\begin{eqnarray}
\label{diag1261}
\xymatrix{
\relax \relax \vee_{\relax \!   n} X  \ar[r]^{\alpha}\ar[d]|{(f)_n }\ar@{}[dr]|{\relax \Downarrow \Gamma ^f _{\alpha}}
& \vee _{\relax \!   m} X \ar[d]|{(f)_m}\\
\relax \vee_{\relax \!   n} Y  \ar[d]|{(g)_n }\ar[r]^{\alpha}\ar@{}[dr]|{\relax \Downarrow \Gamma ^{g}_{\alpha}}
& \vee _{\relax \!   m} Y  \ar[d]|{(g)_m }\\
\relax \vee _{\relax \!   n} Z \ar[r]^{\alpha} & \vee _{\relax \!   m} Z\\
}
\end{eqnarray}
yields $\Gamma^{gf}_{\alpha}$.
\end{theo}

\begin{theo}[\emph{Naturality with respect to tracks in $\trackt$}]
\label{gammaresptrack}
Let $f$ be a map between two weak $\nil_s$ cogroups in the additive track category $\addtrackt$. The unique $\Gamma$-structure associated to $f$ satisfies \emph{naturality with respect to tracks in $\trackt$}: let  $\alpha:\enil^{\nil _s} \lra \fnil^{\nil _s}$, $f:X\lra Y$, $g:Y\lra Z$ in $\trackt$, and $\psi :f\imp g$. Then the pasting in the diagram
\begin{eqnarray}
\label{diag120}
\xymatrix{
\relax \vee _{\relax \!   n} X \ar@/_3pc/[dd]^{\relax \build{\Rightarrow}_{\vee _{\relax \!  n} \psi}^{\boxminus}}_{(g)_n} \ar[r]^{\alpha}\ar[dd]|{(f)_n} \ar[r]^{\alpha}\ar@{}[ddr]|{\relax \Downarrow \Gamma ^f_{\alpha}}
& \vee_{\relax \!   m} X \ar[dd]|{ (f)_m }\ar@/^3pc/[dd]_{\relax \build{\Rightarrow}_{\vee _{\relax \!   m} \psi}^{}}^{(g)_m}\\
&\\
\relax \vee _{\relax \!   n} Y  \ar[r]^{\alpha}& \vee _{\relax \!   m} Y \\
}
\end{eqnarray}
yields $\Gamma_\alpha^{g} $.
\end{theo}
Assume now that all objects are weak $\nil_s$ cogroups in the additive track category $\trackt$, {\it i.e.} we assume that $\addtrackt$ is a  weak $\nil_n$ ringoid. Every map has by the preceding theorem a unique $\Gamma$-structure satisfying  the property $(\Gamma)$. These $\Gamma$-structures have a good behavior with respect to composition in $\trackt$, as we shall see.
Let $\trackt$ be a weak $\nil_s$ ringoid. Theorem \ref{existgamma} builds (according to remark \ref{remremrem}) a correspondence
$$
\Gamma :\trackt \lra  \mathbf{Pseudo} ^\amalg (\nil_s, \trackt )
$$
which to objects associates the chosen weak $\nil_s$ cogroup structure and to each map $f$ in $\trackt$  a coproduct preserving natural transformation $\Gamma ^f$.
Theorem \ref{gammarespmap} says that $\Gamma$ is a functor. Theorem \ref{gammaresptrack} says that $\Gamma$ is actually a $2$-functor. The uniqueness  and existence statements show that $\Gamma$ is in fact an equivalence of $2$-categories (see Theorem \ref{maintheo} and its proof).

\section{\relax Existence and uniqueness of $\Gamma$-structures}
\label{existuniquegamma}

The main point of this section is the proof of Theorem \ref{existgamma}. We proceed in three steps. We first construct an interchange structure, which is a candidate for the $\Gamma$-structure of the theorem. We then show that the constructed interchange tracks are actually $\Gamma$-tracks. The uniqueness is easy, and is derived in the last part of this section.

In this section, we work with some fixed map $f: X\lra Y$ between two weak $\nil_n$ cogroups $X$ and $Y$ in the additive track category $\addtrackt$.

\subsection{Construction of the canonical $\Gamma$-structure of a weak $\nil_n$ ringoid}
\label{gammatrackdefconstruct}

We assume that $\trackt$ has strict coproducts. We will construct an interchange structure associated with $f$ called the \emph{canonical interchange structure} associated with $f$ and denoted by $\Gamma$. We set the interchange tracks $\Gamma^f_\alpha$ associated with projections, inclusions, and fold maps to be identity tracks.

\subsubsection{Additivity $\Gamma$-tracks}
\label{additivity:section}

The diagram
\begin{eqnarray}
\label{diag30}
\xymatrix{\relax
X\ar[r]^{\alpha _{\! n}} \ar[d]_f& \vee _{\! n} X \ar[d]^{(f)_n}\\
Y\ar[r]_{\alpha _{\! n}}&\vee _{\! n} Y \\
}
\end{eqnarray}
is not commutative in $\trackt$ but is commutative  in the homotopy category $\mathrm{ho} (\trackt)$. The map $\alpha _{\! n}$ is described at the beginning of Section \ref{gamma:tracks:section}.
The track set $\mathrm{Track} ((f )_n\alpha _{\! n},\alpha _{\! n} f)$ is therefore non empty. Let $H$ be an element of $\mathrm{Track} ((f)_n\alpha _{\! n}, \alpha _{\! n} f )$. The Abelian group $D(X, \vee _{\! n} Y)$ acts freely and transitively on the set $\mathrm{Track} ((f)_n \alpha _{\! n},\alpha _{\! n} f )$, and the choice of $H$ fixes an isomorphism
\begin{eqnarray}
\sigma _{H}: D(X, \vee _{\! n} Y)\lra \mathrm{Track} ((f)_n \alpha _{\! n}, \alpha _{\! n} f).
\end{eqnarray}
We obtain a diagram
\begin{eqnarray}
\label{diag35}
\xymatrix{\relax &&\\
\relax   X\ar@/^2pc/[rr]^{Id_X}_{\Downarrow \varphi ^{\boxminus} }
\ar@{}[rd]|{\relax \Downarrow H}
\ar[r]^{\alpha _n}\ar[d]_{ f }
& \vee_{\relax \!   n} X  \ar[r]^{r_e} \ar[d]|{(f)_n }& X  \ar[d]^{ f} \\
\relax \ Y
\ar@/_2pc/[rr]_{Id_Y}^{\Downarrow \varphi} \ar[r]_{\alpha _n}& \vee_{\relax \!   n} Y  \ar[r]_{r_e} & Y   &.\\
\relax &&\\
}
\end{eqnarray}
Here the left square commutes because we assume the existence of strict coproducts.
We claim that one can alter $H$ in a unique way by the action of $D(X, \vee _{\! n} Y)$ on the set $\mathrm{Track} ((f)_n \alpha _{\! n}, \alpha _{\! n} f)$, so that pasting in the diagram  (\ref{diag35}) yields the trivial track $f\imp f$, for all $1\leq e \leq n$. The track thus obtained is denoted by $\Gamma _{\! n}^f$
\begin{defi}
The additivity $\Gamma$-track $\Gamma _{\! n}^f$ is the unique track $(f)_n \alpha _{\! n} \imp \alpha _{\! n} f$ that restricts to the trivial track along $r_e:\vee _{\! n} X \lra X$, for all $e\in n $.
\end{defi}
To see that this definition makes sense, we notice that we have a commutative diagram
$$
\xymatrix{
\relax D(X, \vee _{\! n} Y) \ar[rr]^{\sigma _H} \ar[d]|{\Pi _{n}(r_e)_*} & &\relax \mathrm{Track}((f)_n \alpha _{\! n}, \alpha _{\! n} f)\ar[d]|{\Pi _{n}(r_e)_*}\\
\oplus _n D(X,Y)\cong \Pi_ n D(X, Y) \ar[rr]_{~~\relax\prod_{n}\sigma _{(r_e)_* H}}&&\relax \Pi _{n} \mathrm{Track}(f,f)&.\\
}
$$
In this diagram, the top and bottom maps are bijections coming from the structure of a linear track extension. The left vertical map is also an isomorphism, because $\addtrackt$ is an \emph{additive} track extension (thus $D$ is a biadditive bifunctor). It follows that the right vertical map is also an isomorphism. This shows that additivity $\Gamma$-tracks are well defined.

\subsubsection{The negative $\Gamma$-track}

Let $f:X\lra Y$ be a map in $\trackt $. We consider the diagram
\begin{eqnarray}
\label{diag65}
\xymatrix{
\relax X\ar[r]^{\alpha _2}\ar@{}[dr]|{\Downarrow \Gamma_{2}^{f}}\ar[d]_{f}& \vee _2 X \ar[d]|{\relax (f)_2 }\ar[r]^{\relax 1\vee (-1)}& X \ar[d]^{f}\\
Y\ar[r]_{\alpha _2}&\vee _2 Y \ar[r]_{\relax 1\vee (-1)} &Y&.\\
}
\end{eqnarray}
The commutativity of this diagram in the homotopy category forces the existence of  some track $K$ in the following diagram (\ref{diag70}):
\begin{eqnarray}
\label{diag70}
\xymatrix{
\relax X\ar[r]^{\alpha _2}\ar@/^2pc/[rr]_{\relax\Downarrow\varphi ^{\boxminus}}^{0}\ar@{}[dr]|{\Downarrow \Gamma_{2}^{f}}\ar[d]_{f}& \vee _2 X \ar[d]|{\relax  (f)_2 }\ar[r]^{\relax 1\vee (-1)}\ar@{}[dr]|{\Downarrow K}& X \ar[d]^{f}\\
\relax Y\ar@/_2pc/[rr]^{\relax\Downarrow\varphi}_{0}\ar[r]_{\relax \alpha _2}&\vee _2 Y \ar[r]_{\relax 1\vee (-1)} &Y&.\\
&&\\}
\end{eqnarray}
\begin{prop}
\label{uniqgammaneg}
There is a unique track $K$ of the form $0^{\square} \vee \Gamma_{-1}^{f}$ such that the pasting in  diagram (\ref{diag70}) leads to the trivial track $0\imp 0$. This defines and characterizes $\Gamma_{-1}^{f}$. In other words
$$
K \boxbox\Gamma_{2}^{f}=(0^\square , \Gamma_{-1}^{f})\boxbox  \Gamma_{2}^{f}= 0^\square ~~.
$$
\end{prop}
{\it Proof.~} Consider diagram (\ref{diag70}). The category $\trackt$ has strict coproducts and therefore $K=(K_1, K_2)$ with $K_1: f\imp f$ and $K_2: f(-1)\imp (-1)f$. If $K=(K_1, K_2)$ fits in the diagram (\ref{diag70}), then  $K=(0^{\square}, K_2)$ does also. We can thus assume that $K_1= 0^{\square}$.
The condition that the pasting $K \boxbox\Gamma_2^f $ in diagram (\ref{diag70}) yields $0^{\square}:0\imp 0$ determines $\Gamma_2^f \square K$ uniquely. As
$\Gamma_2^f$ is already fixed, this determines $(\alpha _2)^* K$. We consider the commutative diagram
\begin{eqnarray}
\label{diag75}
\xymatrix{
\relax D(X, Y)\times D(X, Y) \ar[r]^{\relax (\alpha _2)^*}\ar[d]^{\relax \sigma_{K_1}\times \sigma_{K_2}}& D(X,Y)\ar[d]^{\relax\sigma_{(\alpha _2)^*}K}\\
\relax \mathrm{Track}(f,f)\times \relax \mathrm{Track}(f(-1),(-1)f)\ar[r]_{~~~\relax (\alpha _2)^*} &\relax ~~~~\mathrm{Track}(f(1,(-1))\alpha _2,(1,(-1))f\alpha _2)&.
}
\end{eqnarray}
The horizontal top map is the sum in the Abelian group $D(X, Y)$. The vertical maps being isomorphisms, we note that if  $\psi_0$ is a fixed element  in $D(X, Y)$ and if we let $\psi$ vary in $D(X, Y)$, then $(\alpha _2)^* (\psi _0, \psi)$ takes all values in $D(X, Y)$ exactly once. Hence $(\alpha _2)^* (\sigma_{K_1}\times \sigma_{K_2}(\psi _0, \psi))$ takes all values in  $\mathrm{Track}(f(1\vee (-1))\alpha _2,(1 \vee (-1))f\alpha _2))$ exactly once. Thus we have proved the existence and uniqueness of $\Gamma_{-1}^{f}:f(-1)\imp (-1)f$ with the prescribed properties.
\fin

\relax
\subsubsection{$\Gamma$-tracks associated to multiplication by a non negative number}

We define a track $\mu _n ^f$ as the result of pasting tracks in the following diagram
\begin{eqnarray}
\label{diag80}
\xymatrix{
\relax
X\ar[r]^{\alpha _n}\ar@{}[dr]|{\relax \Downarrow\ \Gamma_{\! n}^{f}}\ar@/^2pc/[rr]^{\xi _n}_{\Downarrow\varphi ^{\boxminus}}
\ar[d]_{f}& \vee _n  X\ar[r]^{\beta _n}\ar[d]|{\vee_{(f)_n}}&X\ar[d]^{f}\\
\ar@/_2pc/[rr]_{\xi _n}^{\Downarrow\varphi} Y\ar[r]_{\alpha _n}&\relax \vee _{\! n} Y \ar[r]_{\beta _{n}}&Y&.\\
}
\end{eqnarray}

\begin{defi}
We define the track $\Gamma_{\xi _n}^{f}$, more simply denoted by $\mu ^{f}_{n}$, by the equation
\begin{eqnarray}
\Gamma_{\xi _n}^{f}=\mu_{n}^{f}= \Gamma _{\beta _{n}}^{f} \boxbox \Gamma^{f}_{n}
\end{eqnarray}
\end{defi}

\subsubsection{$\Gamma$-tracks associated to multiplication by a negative number}

The negative multiplication $\Gamma$-track $\mu ^f _{\! - \! n}$  is defined as the pasting of tracks in the following diagram:
\begin{eqnarray}
\label{diag85}
\xymatrix{
\relax X\ar[r]^{\xi _n}\ar@/^2pc/[rr]_{\relax \Downarrow \varphi ^{\boxminus}}^{\relax \xi _{-n}}\ar[d]_{f} \ar@{}[rd]|{\Downarrow\mu _n ^f} & X\ar[r]^{\relax \xi_{\! -\! 1}} \ar[d]_{f} \ar@{}[rd]|{\relax \Downarrow \Gamma _{\! -\! 1}^{f}}&X \ar[d]^{f} \\
\relax Y\ar[r]_{\xi _n} \ar@/_2pc/[rr]^{\relax \Downarrow \varphi}_{\relax \xi _{-n}} & Y \ar[r]_{\relax \xi_{\! -\! 1}}  &Y  &.\\
}
\end{eqnarray}
\begin{defi}
We define the track $\Gamma_{\xi _{-n}}^{f}$, simply denoted by $\mu ^{f}_{-n}$, by the equation:
\begin{eqnarray}
\Gamma_{\xi _{-n}}^{f}=\mu_{-n}^{f}= \Gamma _{-1}^{f} \boxbox \mu^{f}_{n}
\end{eqnarray}
\end{defi}

\subsubsection{General $\Gamma$-tracks}
\label{generalgamadef}

We are now ready to define general $\Gamma$-tracks. Let $\alpha:\enil \lra \fnil$, $t:X\lra Y$ be a map in $\trackt$. For each pair $(e,f)$, the
composition $r_{\!f} \alpha  i_e:\z\lra \z$ is the multiplication by a number $\alpha (e, f)$. The $\Gamma$-track $\Gamma_{\alpha}^{f}$ is the unique track $(\vee _{\! m} t) \alpha \imp \alpha (\vee _{\! n} t)$ such that for all $(e, f)\in \underline{n} \times \underline{m}$ the pasting of tracks  in the following diagram yields the track $\mu_{\alpha(e, f)}^{t}$:

\begin{eqnarray}
\label{diag90}
\xymatrix{
\relax X\ar[r]^{i_e}
\ar@/^2pc/[rr]_{\relax \Downarrow  \varphi^{\boxminus}}\ar@/^3pc/[rrr]_{\relax \Downarrow \varphi ^{\boxminus}}^{\xi_{\alpha_{e,f}}}\ar[d]_{t}& \vee _n X \ar@{}[dr]|{\Downarrow K}\ar[d]|{\vee_{\! n} t}\ar[r]^{\alpha}&\relax  \vee _{\! m} X\ar[r]^{r_{\! f}}\ar[d]|{\vee_{\! m} t} & X\ar[d]^{t}\\
\ar@/_3pc/[rrr]^{\relax \Downarrow \varphi}_{\xi_{\alpha_{e,f}}}\ar@/_2pc/[rr]^{\relax \Downarrow \varphi}Y\ar[r]_{i_e}&\vee _n Y \ar[r]_{\alpha}&\relax \vee _{\! m} Y \ar[r]_{\relax r_{\! f}}&Y\\
}
\end{eqnarray}
The uniqueness of such a track follows from two facts. The first one is the assumption that $\trackt$ has strict coproducts, which shows that $\Gamma _{\alpha}^t$ is determined by its various restrictions $i_{\! e}^*\Gamma _{\alpha}^t$, $1\leq e \leq n$. The second point is that  similarly as in the definition of the additivity $\Gamma$-tracks, one shows that  $i_{\! e}^*\Gamma _{\alpha}^t$ is in turn determined by its various corestrictions $(r_f)_* i_{\! e}^*\Gamma _{\alpha}^t$  along the projections $r_f$, for $1\leq f\leq m$ (see also \ref{caracpadretrac}).

By simply checking the definitions:
\begin{prop}
The assignment $\alpha \longmapsto \Gamma_\alpha ^f $ is an interchange structure for $f$.
\end{prop}

\subsection{Existence}

We show in this section that the interchange structure defined in Section \ref{gammatrackdefconstruct} satisfies \emph{property $(\Gamma)$}. This proves the existence part of  Theorem \ref{existgamma}.
\begin{prop}
\label{gammasatis}
Let $\trackt$ be an additive track category with strict coproducts. Let $f: X\lra Y$ be a map of weak $\nil_n$ cogroups. The canonical interchange structure associated to $f$ satisfies the property $(\Gamma)$, and is called the canonical $\Gamma$-structure associated to $f$.
\end{prop}
The proof will be settled through a series of propositions. The first one is: 
\begin{prop}
\label{gamma:one}
For all pair of positive numbers $m, n\in \mathbb{N}$,
$$
\mu ^f _n \boxbox \mu ^f _m = \mu _{nm}^f = \mu _{mn}^f=\mu ^f _m \boxbox\mu ^f _n.
$$
\end{prop}

{\it Proof~.}  We denote by ${\Psi _{n,m}^f}$ the fold map $\vee _m \enil^{\nil_s} \lra \enil^{\nil_s}$.
We have
\begin{eqnarray*}
\mu ^f _n \boxbox \mu ^f _m &=&  \Gamma_{\beta _n}^f\boxbox \Gamma_n^f\boxbox \Gamma_{\beta _m}^f \boxbox\Gamma_m^f\\
&=& \Gamma_{\beta _n}^f\boxbox\Gamma_{\Psi _{n, m}}^f \boxbox (\Gamma _{n}^f)_{i=0}^{m}\boxbox \Gamma_m^f\\
&=& \Gamma_{\beta _{mn}}^f \boxbox \Gamma _{mn}^f\\
&=& \mu _{mn}^f.
\end{eqnarray*}
Here we have used the following two lemmas.
\begin{lem}
\label{lem1}
For all pair of positive numbers $m, n\in \mathbb{N}$, we have
\begin{eqnarray}
\Gamma_n^f\boxbox \Gamma_{\beta _m}^f=\Gamma_{\Psi _{n,m}}^f \boxbox (\Gamma _{n}^f)_{i=0}^{m}.
\end{eqnarray}
\end{lem}

\begin{lem}
\label{lem2}
For all pair of positive numbers $m, n\in \mathbb{N}$, we have
\begin{eqnarray}
\Gamma _{mn}^f=(\Gamma _{n}^f)_{i=0}^{m}\boxbox \Gamma_m^f.
\end{eqnarray}
\end{lem}

\noindent {\it Proof of lemma \ref{lem1} .~} First, one notices that
$\Gamma_{\beta _m}$ is the identity track and therefore:
\begin{eqnarray*}
\Gamma_n^f\boxbox \Gamma_{\beta _m}&=&(\beta _m)^*\Gamma_n^f\\
&=&  (\Psi_{m,n})_*(\Gamma_n^f, \ldots, \Gamma_n^f).
\end{eqnarray*}
Using that $ \Gamma_{\Psi _{n,m}}^{f}=0^\square$, we get
\begin{eqnarray*}
\Gamma_n^f\boxbox \Gamma_{\beta _m}&=& (\Psi_{m,n})_*(\Gamma_n^f, \ldots, \Gamma_n^f) \\
&=& \Gamma_{\Psi _{n,m}}^{f} \boxbox\vee_m \Gamma _{n}^f.
\end{eqnarray*}

\noindent {\it Proof of lemma \ref{lem2} .~}
One only needs to notice that $(\Gamma _{n}^f)_{i=0}^{m} \boxbox \Gamma_m^f$ satisfies the defining property of $\Gamma_{mn}^f$ (see Section \ref{additivity:section}).

\fin

\begin{prop}
\label{gamma:two}
For any non negative number $n\in \mathbb{N}$, we have
\begin{eqnarray}
\mu ^f _{\!-\!1} \boxbox \mu ^f _m =  \mu ^f _m \boxbox \mu ^f _{\!-\!1} .
\end{eqnarray}
\end{prop}
The proof consists of giving a characterization of $\mu ^f _{\!-\!1}\boxbox \mu ^f _m  $ which is also satisfied by $\mu ^f _m \boxbox\mu ^f _{\!-\!1}$. We set $\gamma _{-n}^{f}= \Gamma_{\beta _n}^f\boxbox\vee _{i\leq n} \Gamma_{\! -1}^{f}$. We define $K_n =(\mu_n^f, \gamma_{-n}^{f})$.
\begin{prop}
\label{uniqgammaneg1}
The track $\gamma_{-n}^f$ is the unique track such that the pasting in diagram (\ref{diag1000}) yields the trivial track $0\imp 0$.
\begin{eqnarray}
\label{diag1000}
\xymatrix{
\relax X\ar[r]^{\alpha _{2}}\ar@/^2pc/[rr]_{\relax\Downarrow\varphi ^{\boxminus}}^{0}\ar@{}[dr]|{\Downarrow \Gamma_{2}^{f}}\ar[d]_{f}& \vee _{2} X \ar[d]|{\relax (f)_2 }\ar[r]^{\relax (n \vee -n)}\ar@{}[dr]|{\Downarrow K_n}& X \ar[d]^{f}\\
\relax Y\ar@/_2pc/[rr]^{\relax\Downarrow\varphi}_{0}\ar[r]_{\relax \alpha _{2}}&\vee _2 Y \ar[r]_{\relax (n\vee -n)} &Y\\
&&\\}
\end{eqnarray}
\end{prop}

\proof The proof of this assertion is essentially the same as the proof of Proposition \ref{uniqgammaneg}.
\fin

\begin{lem} We have equalities
$$
\gamma_{-n}^f = \mu ^f _{\!-\!1} \boxbox \mu ^f _m =  \mu ^f _m \boxbox\mu ^f _{\!-\!1}.
$$
\end{lem}

\proof
The proof consists of showing that both  $\Gamma ^f _{\!-\!1} \boxbox \mu ^f _m$ and $ \mu ^f _m \boxbox\mu ^f _{\!-\!1}$ satisfy the characterization of
$\gamma_{-n}^f$ in proposition \ref{uniqgammaneg1}.
We first consider the track $K'_n= (\mu_n^f, \Gamma ^f _{\!-\!1} \boxbox \mu ^f _n )$.
We claim that the pasting $K'_n\boxbox \Gamma_2^f$ in diagram (\ref{diag1010}) yields the trivial track $0\imp0$,
\begin{eqnarray}
\label{diag1010}
\xymatrix{
\relax X\ar[r]^{\alpha _{2}}\ar@/^2pc/[rr]_{\relax\Downarrow\varphi ^{\boxminus}}^{0}\ar@{}[dr]|{\Downarrow \Gamma_{2}^{f}}\ar[d]_{f}& \vee _{2} X \ar[d]|{\relax (f)_2 }\ar[r]^{\relax (n\vee -n)}\ar@{}[dr]|{\Downarrow K'_n}& X \ar[d]^{f}\\
\relax Y\ar@/_2pc/[rr]^{\relax\Downarrow\varphi}_{0}\ar[r]_{\relax \alpha _{2}}&\vee _2 Y \ar[r]_{\relax (n\vee -n)} &Y\\
&&\\}
\end{eqnarray}
Indeed, we have
\begin{eqnarray*}
K'_n\boxbox \Gamma_2^f&=&(\mu_n^f \vee (\Gamma ^f _{\!-\!1} \boxbox \mu ^f _n )) \boxbox \Gamma_2^f  \\
&=& (\Gamma_{\beta_n}^{f}\boxbox \Gamma_{n}^f)\vee (\Gamma ^f _{\!-\!1} \boxbox (\Gamma_{\beta_n}^{f}\boxbox \Gamma_{n}^f) ) \boxbox \Gamma_2^f \\
&=&(\Gamma_{\beta_n}^{f}\vee (\vee_{i\leq n}  \Gamma ^f _{\!-\!1}  ))\boxbox(\Gamma_n^f,\Gamma_n^f) \boxbox \Gamma_2^f \\
&=&(\Gamma_{\beta_n}^{f}\vee (\vee_{i\leq n}  \Gamma ^f _{\!-\!1}  ))\boxbox \Gamma_{2n}^f
\end{eqnarray*}
which is the trivial track $0\imp 0$.
We next consider the track $K_{n}^{''}= (\mu_n^f, \mu ^f _n \boxbox \Gamma ^f _{\!-\!1} )$.
We claim now that the pasting $K_{n}^{''}\boxbox \Gamma_2^f$ in diagram \ref{diag1020} yields the trivial track $0\imp0$,
\begin{eqnarray}
\label{diag1020}
\xymatrix{
\relax X\ar[r]^{\alpha _{2}}\ar@/^2pc/[rr]_{\relax\Downarrow\varphi ^{\boxminus}}^{0}\ar@{}[dr]|{\Downarrow \Gamma_{2}^{f}}\ar[d]_{f}& \vee _{2} X \ar[d]|{\relax  (f)_2 }\ar[r]^{\relax (n\vee -n)}\ar@{}[dr]|{\relax\Downarrow K_n^{''}}& X \ar[d]^{f}\\
\relax Y\ar@/_2pc/[rr]^{\relax\Downarrow\varphi}_{0}\ar[r]_{\relax \alpha _{2}}&\vee _2 Y \ar[r]_{\relax (n\vee -n)} &Y &.\\
&&\\}
\end{eqnarray}
This claim follows from the sequence of equalities
\begin{eqnarray*}
K_{n}^{''}\boxbox \Gamma_2^f&=&(\mu_n^f \vee ( \mu ^f _n \boxbox \Gamma ^f _{\!-\!1} )) \boxbox \Gamma_2^f\\
&=&  (\Gamma_{\beta_n}^{f}\boxbox \Gamma_{n}^f) \vee (\Gamma_{\beta_n}^{f}\boxbox \Gamma_n^f \boxbox\Gamma ^f _{\!-\!1} )\boxbox \Gamma_2^f    \\
&=& (\Gamma_{\beta_{n}}^{f}\vee\Gamma_{\beta_{n}}^{f})\boxbox (\Gamma_n^f , \Gamma_n^f)\boxbox( 0^{\square} ,\Gamma_{\!-\!1}^f)\boxbox \Gamma_2^f
\end{eqnarray*}
but
$$
\Gamma_n^f\Gamma_{\!-\!1}^f= (\Gamma_{\!-\!1}^f)_{i=1}^{n}\Gamma_n^f,
$$
hence
\begin{eqnarray*}
K_{n}^{''}\boxbox \Gamma_2^f&=&(\Gamma_{\beta_{n}}^{f}\vee\Gamma_{\beta_{n}}^{f})\boxbox( (0^{\square})_n ,(\Gamma_{\!-\!1}^f)_n)\boxbox (\Gamma_n^f , \Gamma_n^f)\boxbox \Gamma_2^f\\
&=&\Gamma_{\beta_{n}}^{f}\boxbox ( \Gamma_{\beta_{2}}^{f})_n \boxbox (0^{\square} ,\Gamma_{\!-\!1}^f)_n \boxbox (\Gamma_2^f)_n \boxbox \Gamma_n^f\\
&=&\Gamma_{\beta_{n}}^{f}\boxbox[( \Gamma_{\beta_{2}}^{f}) \boxbox (0^{\square} ,\Gamma_{\!-\!1}^f) \boxbox (\Gamma_2^f)]_n \boxbox \Gamma_n^f,
\end{eqnarray*}
which is the trivial track $0\imp 0$. \fin

Proposition \ref{gamma:one} and \ref{gamma:two} together imply:
\begin{prop}
For all pair of integers $n,m\in \mathbb{Z}$, we have
$$
\mu ^f _{n} \boxbox \mu ^f _m =  \mu ^f _m \boxbox\mu ^f _{n}
$$
\end{prop}

\noindent {\it Proof of proposition \ref{gammasatis}.~} We consider general maps $\alpha, \beta$ in $\nil _s$, with $\alpha :\enil^{\nil _s}\lra \fnil^{\nil _s}$ and $\beta: \fnil^{\nil _s}\lra \gnil^{\nil _s}$.
We assume that $\alpha = \vee_{n} i_e$ where $\alpha _e (e)= \prod_{i} f_i^{n_i} $, with $f_i \in \{1\ldots, m\}$. In the same way, $\beta = \vee_{m} \beta _f$ where $\beta _f (f)= \prod_{j} g_j^{n_j} $, with $g_j \in \{1, \ldots q\}$. Then $\beta\alpha : \enil\lra \gnil$ is such that $\beta\alpha =\vee_{n} (\beta\alpha) _e $
and $(\beta\alpha) _e= \beta \alpha _e$, so that $(\beta\alpha) _e$ maps $e$ to $\prod_{i} (\beta(f_i))^{n_i}$.
We have to show that
$$
\Gamma_\beta^f\boxbox \Gamma_{\alpha _e}^f = \Gamma_{ \beta\alpha_e }^f  ~~.
$$
We begin with the lemma (see \ref{caracpadretrac}):
\begin{lem}
The tracks $\Gamma_\beta^f\boxbox \Gamma_{\alpha _e} ^f$ and $\Gamma_{ \beta \alpha_e }^f$ coincide if and only if the tracks  $\Gamma_{r_g}^f \boxbox \Gamma_\beta^f\boxbox \Gamma_{\alpha _e}^f$ and $\Gamma_{r_g}^f \boxbox \Gamma_{ \beta \alpha_e }^f= \mu_{{(\beta\alpha)}{(e, g)}}$ coincide for all $g\in \{1, \ldots q\}$ (see \ref{generalgamadef} for the definition of ${(\beta\alpha)}{(e, g)}$.
\end{lem}
Next, we  notice  that
\begin{lem}
We have:
$$
\Gamma_{r_g}^f \boxbox \Gamma_{\beta}^f= \vee_m \mu_{n_g}^f
$$
where $n_g= \sum_{g_j=g} n_j$.
\end{lem}
We are thus reduced to show that for all $g\in G$,
$$
\vee_m \mu_{n_g}^f\boxbox\Gamma_{i_e}= \mu_{(\beta\alpha)(e, g)}
$$
But now
\begin{eqnarray*}
\vee_m \mu_{n_g}^f\boxbox\Gamma_{i_e} &=&  \Gamma_{\beta_{n}}^f\boxbox(\mu_{n_g}^f )_{m}\boxbox \Gamma_{i_e}\\
&=& \Gamma_{\beta_{n}}^f \boxbox (\mu_{\beta(f,g)}^f \boxbox \mu_{\alpha_{ef}})_{m}\boxbox \Gamma_{m}^f\\
&=&\mu_{\alpha_{e,f}}
\end{eqnarray*}
and this finishes the proof of Proposition \ref{gammasatis}.\fin

\subsection{Uniqueness}

In this section, we prove:
\begin{theo}
Assume $\trackt$ is an additive track category with strict coproducts. Let $t$ be a map between weak $\nil_n$ cogroups. Two interchange structures associated with $t$ and satisfying property $(\Gamma)$ coincide; there is therefore a unique $\Gamma$-structure associated with $t$.
\end{theo}

The proof consists of a series of lemmas. In fact we will see that under the assumptions, any interchange structure $\tilde{\Gamma}^t$ satisfying property $(\Gamma)$ coincides with  $\Gamma^t$, the canonical one that we have constructed in Section \ref{gammatrackdefconstruct}.
Let $\tilde{\Gamma}^t$ be an interchange structure satisfying property $(\Gamma)$. We begin with the following easy lemma.
\begin{lem}
\label{identitytracktriv}
The $\tilde{\Gamma}$ interchange tracks associated to identities in $\nil _s$ are identities. The $\tilde{\Gamma}$ interchange tracks associated with the trivial map in $\nil_s$ are the trivial tracks.
\end{lem}

\noindent {\it Proof.~}Let $1$ be the identity of $\mathbf{F}_n$. Considering the pasting of $\tilde{\Gamma}$ with itself, we see that
$$
\tilde{\Gamma} _{1}^{t}\boxbox \tilde{\Gamma}_{1}^{t} =\tilde{\Gamma} _{1}^{t}\square \tilde{\Gamma}_{1}^{t} =\tilde{\Gamma}_{1}^{t}
$$
in the group $\mathrm{Track}(t,t)$; therefore $\tilde{\Gamma}_{1}^{t}= 0^{\square}$
\fin

\begin{lem}
The additivity $\Gamma$-tracks are unique, that is $\tilde{\Gamma} _{n}^{t}={\Gamma} _{n}^{t}$ for all $n\geq 0$.
\end{lem}

\noindent {\it Proof.~} Let $n\geq 0$  and $t:X\lra Y$ any map in $\trackt$. The additivity $\Gamma$-track is defined to be the unique one that restricts to the trivial track $t\imp t$ along the retractions $r_e:\vee _{\! n} X \lra X$ for all $e \in n $. Let us see that it coincides with the $\tilde{\Gamma}$-track.
The restriction of the $\tilde{\Gamma}$-track along the maps $r_e:\vee _{\! n} X \lra X$ has to be the $\Gamma$-track associated with the identity, which we have seen (\ref{identitytracktriv}) to be the trivial track. But this property characterizes the additivity $\Gamma$-track, and therefore for all $n$,
$\Gamma_n^{t}= \tilde{\Gamma}_n^t$.
\fin

\begin{lem}
The multiplication $\Gamma$-tracks associated to the multiplication with a positive number is unique.
\end{lem}

\noindent {\it Proof.~} For the multiplication by a positive number, the uniqueness is clear from the diagram (\ref{diag80}) under the hypotheses.
\fin

\begin{lem}
\label{uniqueneg}
The basic negative $\Gamma$-track is unique.
\end{lem}
It follows immediately that:
\begin{cor}
\label{quatretrentequatre}
The $\Gamma$-tracks associated with multiplication maps are unique.
\end{cor}

\noindent {\it Proof of lemma \ref{uniqueneg}.~}
We consider once again the diagram (\ref{diag70})
\begin{eqnarray}
\xymatrix{
\relax X\ar[r]^{\alpha _2}\ar@/^2pc/[rr]_{\relax\Downarrow\varphi ^{\boxminus}}\ar@{}[dr]|{\Downarrow \Gamma_{2}^{t}}\ar[d]_{t}& \vee _2 X \ar[d]|{\relax \vee_{2} t }\ar[r]^{(1,(-1))}\ar@{}[dr]|{\Downarrow K}& X \ar[d]^{t}\\
Y\ar@/_2pc/[rr]^{\relax\Downarrow\varphi}\ar[r]_{\alpha _2}&\vee _2 Y \ar[r]_{(1,(-1))} &Y&.\\
}
\end{eqnarray}
The track $K$ is of the form $(K_1, K_2)$ because $\trackt$ has strict coproducts. By assumption, the \emph{property $(\Gamma)$} is satisfied and this leads
to
\begin{eqnarray}
K\boxbox \Gamma_2 = \Gamma _* ^t = 0^{\square}.
\end{eqnarray}
Because of the definition of a weak $\nil_n$ cogroup, we have
\begin{eqnarray}
(K_1, K_2)=(\Gamma _1 ^t, \Gamma _{\!-1} ^t).
\end{eqnarray}
This fact being granted, we see that  that $\tilde{\Gamma} _{\!-1}$ has to coincide with $\Gamma _{\!-1}$, as these properties determine $\Gamma _{\!-1}$ (see proposition \ref{uniqgammaneg}).
\fin

\begin{prop}
General $\Gamma$-tracks are unique.
\end{prop}

\noindent {\it Proof.~}
According to the remarks following the definition of the general $\Gamma$-tracks (see Section \ref{generalgamadef}), the uniqueness is settled as soon as the multiplication $\Gamma$-tracks are determined, but this precisely the content of corollary \ref{quatretrentequatre}.
\fin

\section{Dual results}
\label{dual}

As we have noticed in \ref{strictproducts}, any additive track category has also a model with strict products. Let $\calc$ be an additive track category with strict products.
One can define in $\calc$ the dual notions of \emph{product preserving pseudofunctor}  $\nil_n ^{op} \sla \calc$ over the algebraic theory of $\nil_s$ groups, \emph{product preserving pseudo natural transformations} between those, and \emph{homotopies of product preserving pseudo natural transformations}. These altogether build a  $2$-category denoted by  $\mathbf{Pseudo}^{\Pi} (\nil_n, \calc)$ and termed the category of $\Pi$-pseudomodels over the algebraic theory $\nil_n^{\mathrm{op}}$.  An object $X$ in $ \calc$, together with a $\Pi$-pseudomodel  $F: \nil_n^{\mathrm{op}} \sla \calc$  such that $F(\mathbb{Z}) = X$, is termed a \emph{weak $\nil_n$ group}.
$\trackt$ is a track category with strict coproducts if and only if its opposite category is a track category with strict products. In this way, we can dualize all our results, and in particular:
\begin{theo}
\label{exist2}
Any object $X$ in an additive track category $\addtrackt$ with strict products admits  canonically the structure of a weak $\nil_n$ group.
\end{theo}
Moreover:
\begin{theo}
\label{uniqnat2}
Let $X, Y$ be two objects in an additive track category $\mathbf{A}$ with strict products, having the weak $\nil_n$ group structures $F_X$ and $F_Y$. Then each map $f:X\lra Y$ extends uniquely to a map
$$
\Gamma ^f: F_X \lra F_ Y
$$
in the category $\mathbf{Pseudo}^{\Pi} (\nil_n , \mathbf{A})$.
\end{theo}
Finally:
\begin{theo}
\label{maintheo2}
The assignment:
\begin{eqnarray*}
T :&\relax \mathbf{A}\lra&\relax \mathbf{Pseudo}^{\Pi} (\nil_n , \mathbf{A}) \\
&X\longmapsto &  F_X
\end{eqnarray*}
extends uniquely to a track functor $T$ which is an equivalence of $2$-categories, the inverse of which is the evaluation of a pseudofunctor $G :\nil_n ^{op}  \sla \mathbf{A}$ on the group $\mathbb{Z}$ for all $n\geq 2$.
\end{theo}

\appendix

\section{Track categories}
\label{trackcat} 
\subsection{Linear track extensions}
\label{lineartrack}

\paragraph{Factorizations}
\label{facto}

Let $\calc$ be a category. The category of factorizations of $\calc$ is the category $\mathcal{F}\calc$ defined by
\begin{itemize}
\item $\mathrm{Ob}(\mathcal{F}\calc)= \mathrm{Mor}(\calc)$

\item For $f, g \in \mathrm{Ob}(\mathcal{F}\calc)= \mathrm{Mor}(\calc)$, a morphism $f\lra g$ is a pair $ (\alpha, \beta)$ fitting in a commutative diagram:
$$
\xymatrix{
{A}{\ar[r]_{\alpha}}&{A'}\\
{B}{\ar[u]^{f}}& {B'}{\ar[l]_{\beta}}{\ar[u]_{g}}\\
}
$$
\end{itemize}

A functor from $D: \mathcal{F}\calc\lra {\mathrm{Ab}}$ is called a natural system on $\calc$.

\begin{ex}
\label{ex:linear:syst}
Let $\calc$ be any category. Any bifunctor $F: \calc ^{\mathrm{op}} \times \calc \lra \abelian$ defines a natural system by
$$
D(X,Y)= F(X, Y), ~ D(\alpha, \beta)= F(\alpha, \beta)
$$
Let $\ab$ be the algebraic theory of Abelian groups (as in Section \ref{algtheo}). Any biadditive bifunctor $F:\ab^{\mathrm{op}} \times \ab$ is of the form
$\mathrm{hom_{\ab}(-, - \otimes M)}$, with $M= D(\mathbb{Z}, \mathbb{Z})$.
\end{ex}

\paragraph{Track categories}
\label{trackcatpar}

A track category $\trackt$ is a category enriched in groupoids. Given two morphisms objects $X, Y\in \trackt$, we have an hom-groupoid $\llbracket X, Y \rrbracket$, its objects are morphisms $f:X\lra Y$ of $\trackt$ and its morphisms $\varphi:f\imp g$ are called tracks from $f$ to $g$. The set of $2$-morphisms from $f$ to $g$ is denoted by $\mathrm{Track}(f,g)$.

The composition of tracks $\eta : f\imp g$ and $\varphi:g\imp h$
is termed \emph{vertical composition} and denoted by $\varphi \square \eta$. The endomorphism-groupoid $\llbracket f, f \rrbracket$ of a $f$  is a group for the vertical composition, termed self tracks of $f$, whose neutral element is denoted by $0^{\square}$. In this group, the inverse of the track $\alpha:f\imp f$ is denoted by $\alpha ^{\boxminus}$.

\begin{ex}
\label{ex:track}
There are numerous examples of track categories, arising in different contexts. Topology provides examples by considering the full subcategories of fibrant-cofibrant objects in \emph{stable} model categories, with tracks being homotopy classes of homotopies. Complete details are presented in \cite{baues1}.
\end{ex}

By definition, the composition
$$
\llbracket A, B \rrbracket  \times \llbracket B, C \rrbracket \lra \llbracket A, C \rrbracket
$$
is a bifunctor. This means that for all $f:A\lra B$ and $g: B\lra C$, functors
$$
g_* : \llbracket A, B \rrbracket \lra \llbracket A, C \rrbracket
$$
and
$$
f^*:\llbracket B, C \rrbracket \lra \llbracket A, C \rrbracket
$$
are defined and commute which each other. In particular, for $\alpha: f_0\imp f_1 \in \llbracket A, B \rrbracket$ and $\beta: g_0\imp g_1 \in \llbracket B, C \rrbracket$ the equation
$$
{g_1}_* \alpha \square {f_0}^* \beta = f_1^* \beta \square {g_0}_* \alpha
$$
holds. This defines the \emph{horizontal composition} of tracks.

From a track category $\trackt$ one can construct two ordinary categories $\trackt _0$ and $\trackt _1$. The category $\trackt _0$ is the underlying category of $\trackt$, obtained by forgetting tracks, while $\trackt _1$ has for objects the morphisms of $\trackt$ and has for morphisms $X\lra Y$ the tracks $f\imp g$ with $f,g:X\lra Y$. The composition
in  $ \trackt _1$ is defined by the horizontal composition of tracks.
The functors \emph{source} and \emph{target} induce two functors
$$\xymatrix{\relax
{ \trackt _1} {\ar@<-0.3ex>[r]}{\ar@<0.3ex>[r]}&{ \trackt _0}~, \\
}
$$
hence we can form an equalizer diagram of categories
$$
\xymatrix{
{ \trackt _1} {\ar@<0.3ex>[r]}{\ar@<-0.3ex>[r]}&{\trackt _0}{\ar[r]}&{\mathrm{ho} (\trackt)}.\\
}
$$
That is, $\mathrm{ho} (\trackt)$ has the same objects as $\trackt$ but its morphisms are obtained from those of $\trackt$ by identifying  morphisms related by a $2$-morphisms in $\trackt$. The category $\mathrm{ho} (\trackt)$ is called the homotopy category of the track category $\trackt$.

A zero object in a category is an object which is both initial and final. All such objects are equivalent. A category with a fixed zero object is called a pointed category.
A strict zero object in a track category $\trackt$ is an object $*$ such that for all object $X$ of $\trackt$, the  hom-groupoids $\llbracket X, * \rrbracket$ and $\llbracket *, X \rrbracket$ are trivial groupoids, with one object and one morphism. The object $*$ is in particular a zero object of the underlying category. A track category with chosen strict zero object is a \emph{pointed track category}. For all objects $X, Y$ in a pointed category, we have a unique  map
$\ast _{X, Y}:X\lra * \lra Y$ with the property that for $g:Y\lra Z$ and $h: W\lra X$
$$
g \ast _{X, Y} = \ast_{X, Z}\mathrm{~and~} \ast _{X, Y}h= \ast _{W, Y}~.$$

\paragraph{Linear track extensions}

A linear track extension $\addtracke$ consists of  the following data:
\begin{itemize}

\item a track category $\mathbf{E}$,

\item a natural system $D$ on $\mathrm{ho} (E)$

\item for all maps $f:X\lra Y$ in $E$, an isomorphism of groups
\begin{eqnarray}
\label{lineartrackpar}
\sigma_f : D_{p(f)}\lra \mathrm{Track}(f,f)
\end{eqnarray}

\item moreover the system of isomorphisms $\{\sigma _f\}$ is required to satisfy:
$$
\forall a\in D_{p(f)}=D_{p(g)}, \forall H \in \mathrm{Track}\{f,g\}, \sigma_f (a) \square H = H \square \sigma _g (a),
$$
$$
\forall \alpha \in D_{p(f)}, g^*\sigma_f (\alpha)= \sigma _{fg} (g^* \alpha),
$$
$$
\forall \beta \in D_{p(f)}, f_*\sigma_g (\beta)= \sigma _{fg} (f_* \beta).
$$
\end{itemize}
The  linear track extension $\addtracke$ is usually depicted by a diagram
$$
\xymatrix{
D+~\ar@{>->}[r]&\relax \mathbf{E}_1\ar@<0.3ex>[r]\ar@<-0.3ex>[r]&\relax\mathbf{E}_0 \ar@{->>}[r]^{p}&\relax \mathrm{ho}(\mathbf{E}).\\
}
$$

\begin{ex}
\label{ex:linar:track}
The example \ref{ex:track} above provides linear track extension with $D(X, Y)= \mathrm{hom}_{\mathrm{ho}(\trackt)} (X, \Omega Y)$.
Complete details are devised in \cite{bauesuniversal}.
\end{ex}

\paragraph{Additive track category}
\label{addittrackcatpar}

\begin{defi}
A linear track extension $\addtracke$
$$
\xymatrix{
D+~\ar@{>->}[r]&\relax \mathbf{E}_1\ar@<0.3ex>[r]\ar@<-0.3ex>[r]&\relax \mathbf{E}_0 \ar@{->>}[r]^{p}&\relax \mathrm{ho}(\mathbf{E})= \calc\\
}
$$
is called additive if
\begin{itemize}
\item the underlying track category $\mathbf{E}$ is pointed,

\item the homotopy category is $\mathrm{ho} (\mathbf{E})$ is additive,

\item $D$ is a biadditive bifunctor.

\end{itemize}
\end{defi}
Note that the zero object of $\mathbf{E}$ will automatically be mapped to a zero object of $\mathrm{ho} (\mathbf{E})$.
\begin{rem}
According to \cite{bauesjib}, an additive track category is essentially determined by its underlying track category. We therefore use the term \emph{additive track category} instead of `additive track extension'.
\end{rem}

\paragraph{Strict products}
\label{strictproducts}

We consider an additive track category $(\trackt, D)$, and we assume \emph{strict} coproducts exist, that is for each pair of objects $(X, Y)$, there is  an object $X\vee Y$ and maps $X\lra X\vee Y, Y \lra X\vee Y$, so that the induced map
$$
 \Psi :\llbracket X\vee Y , Z\rrbracket \lra \llbracket X, Z\rrbracket \times \llbracket Y , Z\rrbracket
$$
is an \emph{isomorphism} of categories for all $Z$. The image of a couple of tracks $(\varphi, \psi) \in \llbracket X, Z\rrbracket _1 \times \llbracket Y , Z\rrbracket _1$ in $\llbracket X\vee Y , Z\rrbracket _1$ by the inverse equivalence  $\Psi ^{-1}$ is denoted by $\varphi \vee \psi$. The fact that we have an equivalence of categories implies the following:
\begin{eqnarray*}
(\varphi \vee \psi) \square (\varphi' \vee \psi')& = &\Psi ^{-1}(\varphi, \psi) \square \Psi ^{-1}(\varphi', \psi')\\
& = &\Psi ^{-1} ((\varphi, \psi) \square (\varphi', \psi'))\\
&= &\Psi ^{-1} (\varphi \square \varphi', \psi \square \psi')\\
&= &(\varphi \square \varphi')\vee (\psi \square \psi').
\end{eqnarray*}
We add for further reference the following easy lemma.
\begin{lem}
\label{caracpadretrac}
Let $E$ be a finite ordered set. A sum of tracks $H= (h_e)_{e\in E}: \vee _{e\in E} f_e \imp \vee _{e\in E} g_e$ is characterized as the unique track that:
\begin{itemize}
\item restricts to $h_e$ along $(r_e)_* (i_e)^*$,
\item restricts to the trivial track $*\imp *$ along $(r_e)_* (i_{e'})^*$ for $e\neq e'$ in $E$.
\end{itemize}
\end{lem}

\paragraph{Cohomology of categories}

Given a small category $\calc$,  we let $N_n (\calc)$ denote the $n^{\mathit{th}}$ stage of the nerve of $\calc$. It consists of all $n$-tuple that form a chain of composable morphisms. Given such a chain $\lambda=(f_1, f_2, \ldots, f_n)\in N_n (\calc)$, we define $\bar{\lambda}$ to be the composition
$$
\bar{\lambda}= f_1 f_2 \ldots f_n \quad .
$$

Let $\calc$ be a small category and let $D$ be a coefficient system on $\calc$. The set of $n$-cochains of $C$ with coefficient in $D$ is the set
$$
\relax C^n (\calc , D)=\left\{\sigma: N_n (\calc) \lra \bigcup _{\lambda \in N_n (\calc)} D(\bar{\lambda}) , \sigma (\lambda) \in D(\bar{\lambda}) \right\}
$$
where $N_n (\calc)$ denotes the nerve of $\calc$.
The object $C^\bullet (\calc , D)$ is actually a graded Abelian group. With the face maps
\begin{eqnarray}
d^i:C^n (\calc , D) \lra C^{n+1} (\calc , D)
\end{eqnarray}
defined by
$$
\left\{\begin{tabular}{lc}$d^0 (\sigma) (f_1, \ldots, f_{n+1}) = (f_1)_*\sigma (f_2,\ldots, f_{n+1})$,&\\
$d^n (\sigma) (f_1, \ldots, f_{n+1}) = (f_{n+1})^*\sigma (f_1,\ldots, f_n)$, &\\
$d^i (\sigma)  (f_1, \ldots, f_{n+1}) = \sigma (f_1,\ldots ,f_i f_{i+1}, \ldots , f_{n+1})$& $\mathrm{if}~ 1<i<n+1,$
\end{tabular}
\right.
$$
and the degeneracies given by
\begin{eqnarray}
s^i (\sigma)(f_1, \ldots, f_{n})= \sigma (f_1, \ldots, f_i, id, f_{i+1}, \ldots f_{n}),
\end{eqnarray}
$C^\bullet( \calc, D)$ becomes a cosimplicial Abelian group. The associated normalized cochain complex is denoted by $NC^\bullet( \calc, D)$ and its homology is by definition the cohomology of $\calc$ with coefficients in $D$:
\begin{eqnarray}
H^* (\calc, D)= H^* (NC^\bullet( \calc, D))~~.
\end{eqnarray}

\paragraph{Classification of linear track categories}

Letting $\calc$ be a small category and $D$ be a natural system over $\calc$, we can define a notion of \emph{maps of linear track extensions} of $\calc$ by $D$. In this way, the linear track extensions  of $\calc$ by $D$ build a category whose connected components is a set $\mathbf{Track} (\calc, D)$ (see \cite{bauesd}).
Moreover,  $\mathbf{Track} (\calc, D)$ is in canonical bijection with $H^3 (\calc, D)$.
From \cite[th. 6.2.1]{bauesjibpira} we have the \emph{strictification theorem},
\begin{theo}
Any linear track extension whose homotopy category has arbitrary (resp. finite) coproducts has in its equivalence class a category with strict (resp. finite) coproducts. A similar  statement holds \emph{mutatis mutandis} by replacing coproducts by \emph{products}. 
\end{theo}

\subsection{Track categories of pseudofunctors}
\label{twocatpseudo}

\subsubsection{Pseudofunctors}
\label{pseudofunct}

The material of this section is adapted from \cite{bauesmuro}. The proofs there can easily be translated to our setting.
Let $\calc$ be and $\calt$ be track categories. A pseudofunctor $\calc \sla \calt$ is an assignment of objects, maps, and tracks together with additional tracks
$$
\varphi_{f,g} : \varphi (f) \varphi (g) \imp \varphi (fg) \mathrm{~and~} \varphi_{X}: \varphi (1_{X}) \imp 1_{\varphi _X}
$$
for all objects $X$ and composable maps $\bullet \stackrel{g}{\lra} \bullet \stackrel{f}{\lra} \bullet$. These tracks must satisfy the following conditions.
For all composable $\bullet \stackrel{g}{\lra} \bullet \stackrel{f}{\lra} \bullet \stackrel{h}{\lra} \bullet$
\begin{itemize}
\item the pasting in the diagrams  (\ref{diag2pseudo}) and (\ref{diag444pseudo}) is the identity track,
\item for any composable tracks $\alpha$,  $\beta$  in $\calc$, the composition in the diagram (\ref{diag4pseudo}) is the track $\varphi(\alpha \beta)$,
\item $\varphi$ preserves vertical composition of tracks,
\item $\varphi$ preserves identity tracks.
\end{itemize}
We leave it to the reader to write the formal equations corresponding to these conditions.
\begin{eqnarray}
\label{diag2pseudo}
\xymatrix{
\relax &&\relax &\\
\relax\varphi (X)\ar@/^1.5pc/[r]^{1_{\relax\varphi (X)}}_{\relax \Downarrow \varphi _X ^{\boxminus} }\ar[r]_{\varphi (1_X)}
\ar@/_1.5pc/[rr]_{\varphi (f)}^{\Downarrow \varphi_{f, 1_X}} &\relax \varphi (X)\ar[r]^{\varphi (f)}&\relax\varphi (Y)&\relax \varphi (X)\ar[r]^{\varphi (f)}\ar@/_1.5pc/[rr]_{\varphi (f)}^{\Downarrow \varphi _{1_Y, f}} &\relax \varphi (Y)\ar@/^1.5pc/[r]^{1_{\relax\varphi (Y)}}_{\varphi _Y ^{\Downarrow\boxminus}}\ar[r]&\relax\varphi (Y)\\
\relax &\relax &\\
}
\end{eqnarray}

\begin{eqnarray}\label{diag4pseudo}
\xymatrix{\relax
\relax\bullet\ar@{}[d]|{\relax\build{\Rightarrow}_{\beta}^{}}\ar@/^1pc/[d]^{\varphi (g')}\ar@/_4pc/[dd]_{\varphi (fg)}^{\relax\build{\Rightarrow} _{\varphi_{f,g}^{\boxminus}}^{}}\ar@/^4pc/[dd]^{\varphi (f'g')}_{\relax\build{\Rightarrow}_{\varphi_{f',g'}}^{}}\ar@/_1pc/[d]_{\varphi (g)}\\
\relax\bullet\ar@{}[d]|{\relax\build{\Rightarrow}_{\alpha}^{}}\ar@/^1pc/[d]^{\varphi (f')}\ar@/_1pc/[d]_{\varphi (f)} &\\
\bullet}
\end{eqnarray}

\begin{eqnarray}
\label{diag444pseudo}
\xymatrix{\relax
&&\bullet\ar[d]|{\varphi(h)}\ar@/_2pc/[dd]_{\varphi(gh)}^{\relax\build{\varphi_{g, h}^\boxminus}_{\Rightarrow}^{}}\ar`r[rr]`d[rrddd]^{\varphi(fgh)}_{\relax\build{\Rightarrow}_{\varphi_{fgh}}^{}}[ddd]\ar`l^d[ll]`[ddd]_{\varphi(fgh)}^{\relax\build{\varphi_{f,gh}^\boxminus}_{\Rightarrow}^{}}[ddd] &&\\
&&\bullet\ar@/^2pc/[dd]^{\varphi(fg)}_{\relax\build{\varphi_{f, g}}_{\Rightarrow}^{}}\ar[d]|{\varphi(g)} &&\\
&&\bullet\ar[d]|{\varphi(f)}&&\\
&&\bullet&&\\
}
\end{eqnarray}
\subsubsection{Pseudo-natural transformations}

Given two pseudofunctors
$$
\varphi, \psi:\calc \sla \calt~~,
$$
a pseudo natural transformation $\alpha: \varphi \lra \psi$ is a collection of maps
$$
\alpha_X : \varphi (X) \lra \psi (X)
$$
for all $X$ in $\calc$
and  for all maps $f: X\lra Y$ in $\calc$, a collection of tracks
$$
\alpha _f: \alpha_Y \varphi (f) \imp \psi (f)\alpha_X,
$$
such that
\begin{itemize}
\item for all $f,g:X\lra Y$ in $\calc$ and all tracks $\gamma : f\imp g$, pasting in diagram (\ref{diag5pseudo}) yields $\alpha_g$,
\item for any composable maps
$$
X \stackrel{f}{\lra} Y \stackrel{g}{\lra} Z
$$
in  $\calc$, pasting in the diagram (\ref{diag6pseudo}) yields $\alpha_{fg}$,
\item for all $X$ in $\calc$, pasting in the diagram (\ref{diag7pseudo}) yields the identity track $0^{\square}_{\alpha _X}$.
\end{itemize}
Let $\alpha: \varphi \lra \psi$ and $\beta : \psi \lra \xi$ be two pseudo natural transformations. We define the composite pseudo natural transformation
$$
\beta \alpha : \varphi \lra \xi
$$
as the assignment
$$
(\beta \alpha)_X =\beta _X \alpha _X : \varphi (X) \lra \xi (X)
$$
for all $X$ in $\calc$,
and for all maps $f: X\lra Y$ in $\calc$ the tracks
$$
(\beta\alpha )_f: \beta _Y\alpha_Y \varphi (f) \imp \psi (f) \beta _X\alpha_X
$$
are given by pasting in diagram  (\ref{diag8pseudo}).

The fact that $\beta \alpha$ is again a pseudo natural transformation is straightforward, as well as the fact that this composition is associative. Moreover, for all pseudofunctor $\varphi:\calc \sla \calt$, there is  an identity pseudo natural transformation  $1_\varphi : \varphi  \lra \varphi$ for the composition of pseudo natural transformations, for which all maps $(1_\varphi )_X : \varphi (X)\lra \varphi (X)$ are identities, and all tracks $(1_\varphi) _f: \varphi (f)\imp \varphi (f)$ are identity tracks.

\begin{table}[h]
\begin{eqnarray}
\label{diag8pseudo}
\xymatrix{
\relax \relax \varphi (X)  \ar[r]^{\varphi (f)} \ar[d]_{\alpha _X }\ar@{}[dr]|{\relax \Downarrow \alpha _f}
& \varphi (Y) \ar[d]^{\alpha _Y}\\
\psi (X) \ar[d]_{\beta _X }\ar[r]^{\psi (f)}\ar@{}[dr]|{\Downarrow\beta _f}
& \xi (X)  \ar[d]^{\beta _Y}\\
\xi (x)\ar[r]^{\xi (f)} & \xi (Y)\\
}
\end{eqnarray}
\end{table}

\begin{prop}
Let $\calc$ be a small track category and $\calt$ be any track category. Then the pseudofunctors $\calc\sla \calt$ and their pseudo natural transformations build a category denoted by $\mathbf{Pseudo} (\calc, \calt)$.
\end{prop}

\begin{eqnarray}
\label{diag5pseudo}
\xymatrix{\relax
\ar@{}[drr]|{\relax \Downarrow {\varphi (\gamma) ^{\boxminus}}}
\relax& &\\
\relax \varphi (X)\ar@/^2pc/[rr]^{\varphi(g)} \ar[rr]|{\varphi (f)}\ar[d]_{\alpha _X} \ar@{}[drr]|{\relax \Downarrow \alpha _f}
& &\varphi (Y) \ar[d]^{\alpha _{Y} }\\
\psi (X)\ar@/_2pc/[rr]_{\psi(g)}\ar[rr]|{\psi (f)}&& \psi (Y) \\
\ar@{}[urr]|{\relax \Downarrow {\varphi (\gamma)}} &&\\
}
\end{eqnarray}

\begin{eqnarray}
\label{diag6pseudo}
\xymatrix{
\relax &&\\
\relax \varphi (X)\ar@/^2pc/[rr]^{\varphi(gf)} \ar[r]^{\varphi (f)}\ar[d]_{\alpha _X}\ar@{}[rru]|{\Downarrow \varphi_{g,f}^{\boxminus}}
\ar@{}[dr]|{\relax  \Downarrow \alpha _f }
& \varphi (X)  \ar[r]^{\varphi (g)} \ar[d]|{\alpha _Y}
\ar@{}[dr]|{\relax \Downarrow \alpha _g}
&\varphi (Z)  \ar[d]^{\alpha _Z} \\
\psi (X)  \ar@{}[rrd]|{\Downarrow \psi_{g,f}}\ar@/_2pc/[rr]_{\psi (gf)}\ar[r]^{\psi (f)}& \psi (Y)  \ar[r]^{{\psi (g)}} &\psi (Z)   \\
\relax &&\\
}
\end{eqnarray}

\begin{eqnarray}
\label{diag7pseudo}
\xymatrix{
\ar@{}[drr]|{\relax \Downarrow {\varphi _X} ^{\boxminus}}
\relax& &\\
\relax \varphi (X)\ar@/^2pc/[rr]^{1_\varphi(X)} \ar[rr]|{\varphi (1_X)}\ar[d]_{\alpha _X} \ar@{}[drr]|{\relax \Downarrow \alpha _{1_X}}
& &\varphi (X) \ar[d]^{\alpha _{X} }\\
\psi (X)\ar@/_2pc/[rr]_{1_\psi(X)}\ar[rr]|{\psi (1_X)}&& \psi (X) \\
\ar@{}[urr]|{\relax \Downarrow {\psi _X}} &&\\
}
\end{eqnarray}

\subsubsection{The track category of pseudofunctors}

Let  $\calc$ be a small track category and $\calt$ be any track category. Let $\alpha, \beta : \varphi \lra \psi$ be a morphism in $\mathbf{Pseudo} (\calc, \calt)$. A track $H: \alpha \imp \beta$ is a collection of tracks
$$
H_X : \alpha_X \imp \beta _X
$$
for all $X\in \calc$, such that
\begin{itemize}
\item for all $f: X \lra Y$ in $\calc$, pasting in  diagram (\ref{diag9pseudo}) yields $\beta_f$,
\item for all $X$ in $\calc$, pasting in  diagram (\ref{diag10pseudo}) yields $\beta _{1_X}$,
\end{itemize}
\begin{eqnarray}
\label{diag9pseudo}
\xymatrix{
\relax \varphi (X) \ar@/_2pc/[d]^{\relax \build{\Rightarrow}_{H_X^{\boxminus}}^{}}_{\beta _X} \ar[r]^{\varphi (f)}\ar[d]|{\alpha _X} \ar[r]^{\varphi (f)}\ar@{}[dr]|{\relax \Downarrow \alpha _f}
& \varphi (Y) \ar[d]|{\alpha _Y }\ar@/^2pc/[d]_{\build{\Rightarrow}_{H_Y}^{}}^{\beta _Y}\\
\relax \psi (X)  \ar[r]^{\psi (f)}& \psi (Y)&. \\
}
\end{eqnarray}

\begin{eqnarray}
\label{diag10pseudo}
\xymatrix{
\relax \varphi (X)\ar@/_2pc/[d]^{\relax \build{\Rightarrow}_{H_X^{\boxminus}}^{}}_{\beta _X} \ar[rr]^{\varphi (1_X)}\ar[d]|{\alpha _X} \ar@{}[drr]|{\relax \Downarrow \alpha _{1_X}}
& &\varphi (X) \ar[d]|{\alpha _{X} }\ar@/^2pc/[d]_{\build{\Rightarrow}_{H_Y}^{}}^{\beta _Y}\\
\psi (X)\ar[rr]_{\psi (1_X)}&& \psi (X) \\
\\
}
\end{eqnarray}
The tracks of pseudo natural transformations $\mathbf{Pseudo} (\calc, \calt)$ are actually the $2$-morphisms of a track structure on $\mathbf{Pseudo} (\calc, \calt)$.
\begin{prop}
Let $\calc$ be a small track category and $\calt$ be any track category. Then the pseudofunctors $\calc\sla \calt$, their pseudo natural transformations, and tracks of pseudo natural transformations give $\mathbf{Pseudo} (\calc, \calt)$ the structure of a track category.
\end{prop}

\subsection{Basic results on pseudofunctors (after \cite{bauesmuro})}
\label{pseudopseudo}

Let $\calc$ and $\calt$ be two track categories and $\varphi: \calc \sla \calt$ a pseudofunctor.

The pseudofunctor $\varphi$ is reduced if
$$
\varphi (id_X) = id_{\varphi (X)} \mathrm{~and ~} \varphi  (0_ X^\square) = 0_{\varphi (X)}^\square
$$
for all objects $X$ in $\calc$.

We now assume that $\calc$ and $\calt$ have a strict zero object. The pseudofunctor $\varphi$ is normalized at zero maps if
\begin{itemize}
\item $\varphi(\ast ) =\ast $,
\item $\varphi (\ast _{X, Y}) = \ast _{\varphi(X), \varphi (Y)}$,
\item $\varphi_{f, \ast }=\varphi _{\ast , f}$ is the trivial track of the zero map.
\end{itemize}

Pseudofunctor reduced and normalized at zero objects are called completely reduced.

\begin{prop}
Let $\varphi: \calc \sla \calt$ be a pseudofunctor. We consider a collection of tracks $\xi=\{\xi _f: \varphi(f)\imp \varphi^{f}\}_{f\in \mathrm{Mor} ~{\calc}}$.
The correspondence $\varphi^\xi $ that associates to objects $X$ in $\calc$
$$
X \longmapsto \varphi^\xi (X)= \varphi (X)\quad ,
$$
to maps $f: X\lra Y$ in $\calc$
$$
f \longmapsto \varphi^\xi (f)= \varphi ^f,
$$
to tracks $\alpha:f\imp g$
$$
\varphi^\xi (\alpha)= \xi _g \varphi (\alpha) (\xi_f)^\boxminus
$$
together with the tracks
$$
\varphi^\xi _{f, g}= \xi _{fg}\square \varphi_{f, g} \square (\xi _f^\boxminus,\xi _g^\boxminus ),
$$
and
$$
\varphi^\xi _{X}= \varphi_{X} \xi _{1_X}^\boxminus
$$
define
\begin{itemize}
\item a pseudofunctor $\varphi^\xi $, and
\item a pseudo natural transformation $t_\xi: \varphi \lra \varphi ^\xi$ that is the identity of objects.
\end{itemize}
Moreover,
\begin{itemize}
\item given a second collection $\xi '=\{\xi' _f: \varphi^f\imp (\varphi')^{f}\}_{f\in \mathrm{Mor} ~{\calc}}$, we may consider the  collection
$\xi '*\xi=\{(\xi*\xi ') _f: \varphi^f\imp (\varphi)^{f}\imp (\varphi')^f \}_{f\in \mathrm{Mor} ~{\calc}}$. W)e have
$$ t_{\xi'} t_\xi = t_{\xi'*\xi} .$$
\item if $\calc$ and $\calt$ have strict coproducts (resp. products), and $\varphi: \calc\lra \calt$ is coproduct (resp. product) preserving, then $\varphi^\xi$ is again  coproduct  (resp. product) preserving and  $t_\xi$ is a coproduct (resp. product)  preserving pseudo natural transformation.
\end{itemize}
\end{prop}

The proof of this proposition is straightforward, and yields easily the following corollaries.
\begin{cor}
Any pseudofunctor  $\varphi: \calc \lra \calt$ is naturally homotopic to a reduced pseudofunctor. If $\calc$ and $\calt$ have strict coproducts (resp. products), and $\varphi: \calc\lra \calt$ is coproduct (resp. product) preserving, then $\varphi^\xi$ is naturally homotopic to a reduced coproduct  (resp. product) preserving pseudofunctor through a coproduct (resp. product) preserving pseudo natural transformation.
\end{cor}

\begin{cor}
Any pseudofunctor  such that $\psi (*) =*$ is naturally isomorphic a completely reduced pseudofunctor. If $\calc$ and $\calt$ have strict coproducts (resp. products), and $\varphi: \calc\lra \calt$ is coproduct (resp. product) preserving, then $\varphi^\xi$ is naturally homotopic to a completely reduced coproduct  (resp. product) preserving pseudofunctor through a coproduct (resp. product) preserving pseudo natural transformation.
\end{cor}

\end{document}